\documentclass[a4paper, 10pt]{article}

\usepackage{graphicx}
\usepackage{multirow}
\usepackage{amsmath,amssymb,amsfonts}
\usepackage{amsthm}%
\usepackage{mathrsfs}%
\usepackage[title]{appendix}%
\usepackage{xcolor}%
\usepackage{textcomp}%
\usepackage{manyfoot}%
\usepackage{booktabs}%
\usepackage{algorithm}%
\usepackage{algorithmicx}%
\usepackage{algpseudocode}%
\usepackage{listings}%
\usepackage{url}
\usepackage{relsize}
\usepackage{enumerate}
\usepackage{hyperref}
\usepackage{color}
\usepackage{lipsum}
\usepackage{float}
\usepackage{xpatch}
\usepackage{xifthen}
\usepackage{verbatim}
\usepackage{todonotes}
\usepackage{bbm}
\usepackage{epsfig}
\usepackage{authblk}
\usepackage{enumitem}			
\usepackage{bbold}
\usepackage{subfigure}
\usepackage{cleveref}
\usepackage{pgfplots} 				
\usepackage{tikz}
\usepackage{tikzscale}
\usepackage{environ}

\newcommand{\R}{\mathbb{R}}
\newcommand{\N}{\mathbb{N}}

\newcommand{\inner}[2]{\ifthenelse{\equal{#2}{}}{\left\langle\cdot,\cdot\right\rangle_{#1}}{\left\langle#2\right\rangle_{#1}}}
\newcommand{\scalprod}[2]{ \left< #1 , #2 \right> }
\newcommand{\norm}[2]{\ifthenelse{\equal{#2}{}}{\left\|\cdot\right\|_{#1}}{\left\|#2\right\|_{#1}}}
\newcommand{\seminorm}[2]{\ifthenelse{\equal{#2}{}}{\left|\cdot\right|_{#1}}{\left|#2\right|_{#1}}}

\newcommand{\mynorm}[1]{\left\|#1 \right\|}

\newcommand{\nsbar}{\mathcal H_k (\bar \Omega)}
\newcommand{\calh}{\mathcal H}
\newcommand{\cO}{\mathcal O}
\newcommand{\cP}{\mathcal P}

\newcommand{\bx}{{\bf x}}
\newcommand{\bOmega}{\boldsymbol{\Omega}}

\DeclareMathOperator*{\myspan}{span}

\newtheorem{theorem}{Theorem}
\newtheorem{prop}[theorem]{Proposition}
\newtheorem{cor}[theorem]{Corollary}

\newtheorem{rem}[theorem]{Remark}

\usepackage[normalem]{ulem}			

\graphicspath{{./Figures/}}
\pgfplotsset{compat=1.11}
\newlength\fwidth

\theoremstyle{thmstyleone}%
\theoremstyle{thmstyletwo}%

\theoremstyle{thmstylethree}%

\newcommand{\hide}[1]{ }

\newenvironment{delayedproof}[1]
 {\begin{proof}[\raisedtarget{#1}Proof of \Cref{#1}]}
 {\end{proof}}
\newcommand{\raisedtarget}[1]{%
  \raisebox{\fontcharht\font`P}[0pt][0pt]{\hypertarget{#1}{}}%
}

\begin{document}

\title{Kernel-based Greedy Approximation of Parametric Elliptic Boundary Value Problems}

\author[1]{Bernard Haasdonk \thanks{haasdonk@mathematik.uni-stuttgart.de}}
\author[2]{Tizian Wenzel \thanks{wenzel@math.lmu.de}}
\author[3]{Gabriele Santin \thanks{gabriele.santin@unive.it}}

\affil[1]{Institute of Applied Analysis and Numerical Simulation, University of Stuttgart, Germany}
\affil[2]{Department of Mathematics, LMU Munich, Germany}
\affil[3]{Department of Environmental Sciences, Informatics and Statistics, Ca' Foscari University of Venice, Italy}

\maketitle 

\begin{abstract}
  We recently introduced a scale of kernel-based greedy
  schemes for approximating the solutions of elliptic boundary value problems.
  The procedure is based on a generalized interpolation
  framework in reproducing kernel Hilbert spaces
  and was coined PDE-$\beta$-greedy procedure, where the
  parameter $\beta\geq 0 $ is used in a greedy selection criterion and
  steers the degree of function adaptivity.
  Algebraic convergence rates have been obtained for Sobolev-space
  kernels and solutions of finite smoothness.
  We now report a result of exponential convergence rates for the case of
  infinitely smooth kernels and solutions.
  We furthermore extend the approximation scheme to the case of parametric PDEs by the use of
  state-parameter product kernels. 
  In the surrogate modelling context,  
  the resulting approach can be interpreted as an {\it a priori}
  model reduction approach, as no solution snapshots need to be precomputed. 
  Numerical results show the efficiency of the approximation
  procedure for problems which occur as challenges for
  other parametric MOR procedures:
  non-affine geometry parametrizations, moving sources or
  high-dimensional domains.\\

\end{abstract}

\section{Introduction}\label{sec:intro}

We consider general linear partial differential equation (PDE) boundary value problems (BVPs) in the strong form
on arbitrary bounded Lipschitz domains.
Apart from classical grid-based discretization techniques, collocation
approaches have frequently been discussed in the literature. 
In kernel-based approximation of PDEs, a well known approach is Kansa's method \cite{Kansa1990a,Schaback2007}, 
which results in a non-symmetric system matrix.
In contrast, changing the ansatz-space for the approximants
leads to symmetric collocation procedures
\cite{fasshauer1997} which have computational and analytical advantages.
Especially these symmetric approaches can be cast in the setting of
generalized interpolation \cite[Chap.\ 16]{wendland2005scattered},
which is an elegant and powerful framework covering many optimal
recovery problems.
In the context of approximating BVPs by collocation, the
problem of positioning the collocation points is a central
aspect determining the quality of the approximation. 
We refer to the PDE-$P$-greedy procedure \cite{schaback2019greedy},
which provides a data-independent recipe for placing the points. 
As a generalization, we recently
presented a scale of meshless greedy kernel-based collocation techniques
\cite{PDEgreedy2025} denoted as PDE-$\beta$-greedy procedures. 
The approximation spaces are incrementally constructed by
carefully collecting Riesz-representers of (derivative operator)
point-evaluation functionals.
The scale of methods is parametrized by a parameter $\beta\in [0,\infty]$ and
naturally generalizes existing approaches of
PDE approximation such as the PDE-$P$-greedy 
as well as function approximation techniques such as the $P$-greedy 
\cite{marchi2005optimal},
$f$-greedy \cite{mueller2009komplexitaet,SchWen2000},
$f/P$-greedy \cite{mueller2009komplexitaet,Wirtz2013},
$f \cdot P$ (also denoted power scaled residual (psr) method) \cite{Dutta2020},
or the scale of  $\beta$-greedy approaches \cite{wenzel2022analysis}.
We remark that point-selection methods are widely studied for the Kansa collocation method, see~\cite{Ling2004,Ling2006,Ling2008,Ling2009,Lee2009}.
Moreover, another popular technique for the approximation of PDE solutions by kernels is the RBF-FD method, which is based on local differentiation 
stencils similarly to classical finite difference schemes. For this algorithm, the selection of interpolation and collocation points is of central importance, 
and several techniques have been introduced~\cite{Vlasiuk2018,Davydov2011,Suchde2023,Oanh2017,Lehto2017,Shankar2018}.

General 
convergence statements for greedy generalized interpolation methods
have been derived \cite{Albrecht2024}. 
In contrast to those non-quantitative statements, 
assuming well-posedness and a stability estimate of the given
BVP, we can rigorously prove convergence rates of the
resulting approximation schemes \cite{PDEgreedy2025}.
Interestingly,
those rates show that it is possible to break the curse of
dimensionality and potentially reach high input dimensions.
In the (non-PDE) function approximation context, it could be shown that
those rates (for the $f$-greedy procedure) are optimal for certain
kernels and target functions \cite{santin2024ontheoptim}.
For domains in small input space dimensions the schemes
can be experimentally compared with, e.g., standard finite element methods (FEM).
The strength of the procedure, however, is the ease of treating
high-dimensional input space dimensions due to the mesh-independence and
the omission of spatial integrals.
In the current article,
we 
report a result of exponential convergence rates
for the procedure in the case of smooth solutions and infinitely smooth kernels.
We then extend the scheme in order to approximate
parametric BVPs.
The essential point here is the use of
combined parameter-state kernel functions.
Approximation of parametric problems is the
central task in parametric model order
reduction (MOR) \cite{BGW15} and
reduced basis (RB) methods \cite{Haasdonk2017,Hesthaven2016}.
The latter especially represent so
called {\em a posteriori} MOR procedures, 
that make use of an
offline-online decomposition, where in the offline-phase snapshots of
solutions are computed in order to construct a parametric surrogate, that can be
efficiently queried in an online-phase.
For the case of collocation discretizations we want to mention the 
Reduced Collocation approaches \cite{CG13} and emphasize that the online-phase
for such models is non-trivial, as it requires reduced
system assembly (possibly involving additional
interpolation steps), and solving of corresponding (least-squares) systems.
In contrast, {\em a priori} MOR approaches aim at providing a global parametric surrogate as
a function, that can be trivially sampled. This eliminates the classical offline-online separation, as the
online-phase is essentially trivialized, and the full model computation is done in the
offline phase.  
In this sense our overall procedure can be interpreted as an {\em a priori}
surrogate modelling approach.
Approaches such as Proper Generalized Decomposition (PGD) \cite{nouy2010}
or Physics Informed Neural Networks (PINNS) \cite{PINNref} are
alternative methods that can provide global surrogates for parametric
problems. 
For fixed (non-greedy) collocation point sets,
recently a constrained minimum norm approximation strategy using Gaussian Processes (GPs) has
been suggested and analyzed for nonlinear PDEs \cite{Chen2021,OwhadiErrorAnalysis}.
For linear PDEs the underlying algebraic linear equation system structure is
identical to ours.
In this viewpoint, our method can be interpreted as a
way of providing problem-adapted collocation point sets
for the GP-approach.
With our approach, we can especially address nontrivial
parametric geometries even allowing topology changes,
non-separable parameter dependencies such as moving sources, or high-dimensional
domains which are inaccessible for standard grid-based methods such as FEM.
These problem types pose challenges or obstacles for standard linear-subspace model
reduction techniques such as RB methods.

The structure of the article is as follows.
We start with recalling the
kernel-based (non-parametric) PDE-greedy schemes in Sec.\ \ref{sec:PDEgreedy}.
We review existing algebraic convergence statements and present an 
exponential convergence 
result in Sec.\ \ref{sec:convergence-rates}.
In Sec.\ \ref{sec:muPDEgreedy} we extend the scheme in order to approximate
parametric PDEs.
Numerical results
in Sec.\ \ref{sec:experiments}
show the benefits and limitations of the approach.
We conclude with some
opportunities for future work in Sec.\ \ref{sec:conclusion}.
An appendix contains the proof details in order
to keep the main text concise.

\section{PDE-Greedy Kernel Methods}

\label{sec:PDEgreedy}

\subsection{PDEs}
We consider BVPs for
linear elliptic second order differential equations in the strong form: For a
given bounded (open) domain $\Omega \subset \R^d$ with boundary $\partial \Omega$,
right-hand side (RHS) source function $f\in C^0(\Omega)\cap L^\infty(\Omega)$ and boundary
value function $g\in C^0(\partial \Omega)$, we want to find a solution
$u\in C^2(\Omega)\cap C^0(\bar \Omega)$
s.t.
$$
Lu(x) = f(x),\quad x \in \Omega, \qquad Bu(x) = g(x),\quad x \in \partial\Omega.
$$
Here  $L:C^2(\Omega) \rightarrow C^{0}(\Omega)$ is assumed to be a differential
operator with coefficients in $C^\infty(\Omega)\cap L^\infty(\Omega)$
and $B:C^1(\Omega) \rightarrow C^0(\partial \Omega)$ is a suitable boundary
value operator, e.g.\ extracting Dirichlet, Neumann or Robin boundary values, possibly with different conditions on different parts of the boundary.
We assume the BVP to be well-posed, i.e.\ for any suitably regular $f$ and $g$ there exists
a unique classical solution $u$, which depends continuously on the data
via
\begin{equation}
\mynorm{u}_{L^\infty(\bar \Omega)} \leq C_L \mynorm{f}_{{L^\infty(\Omega)}} + C_B
\mynorm{g}_{L^\infty(\partial \Omega)},
\label{eqn:well-posedness-estimate}
\end{equation}
with $C_L,C_B$ independent of $u$. 
For uniformly elliptic problems this is typically guaranteed by a suitable maximum
principle \cite{renardy2006introduction,PDEgreedy2025}.

\subsection{Kernels and RKHSs}

The approximation framework is formulated in reproducing kernel Hilbert spaces (RKHSs).
For this we consider a kernel, i.e.\ a symmetric
function $k:\bar \Omega \times \bar \Omega\rightarrow \R$,
which is assumed to be strictly positive definite meaning that
all kernel matrices for pairwise distinct point sets
are assumed to be positive definite.
Such a kernel induces an RKHS denoted $\nsbar \subset \{h\, |\, h:\bar \Omega \rightarrow \R \}$
with the property that $k(x,\cdot)\in \nsbar$ 
and $\scalprod{h}{k(x,\cdot)}_{\nsbar} = h(x)$ for
all $x\in\bar \Omega$ and $h\in\nsbar$ (reproducing property).
For any continuous linear functional $\lambda \in \nsbar'$ we denote
$v_\lambda\in \nsbar$ as its Riesz-representer.
Let  $\delta_x$ denote the
point evaluation functional in $x\in \bar \Omega$,
i.e.\ $\delta_x(h) = h(x)$ for all $h\in \nsbar$.
The definition of an RKHS implies that the 
point evaluation functionals are continuous and
$v_{\delta_x} = k(x,\cdot)$ is the corresponding Riesz-representer.
If $k$ is chosen as a Mat\'ern
or Wendland kernel,
the corresponding RKHS is norm-equivalent to a certain Sobolev space
\cite{wendland2005scattered}.
Especially, we will apply the quadratic Mat\'ern kernel in the experiments, i.e.
\begin{equation}
k^M_\varepsilon(x,x') := (3+3 \varepsilon \mynorm{x-x'}+\varepsilon^2\mynorm{x-x'}^2)
\exp(-\varepsilon \mynorm{x-x'}), \;x,x'\in\bar \Omega,\label{eqn:matern}
\end{equation}
with shape parameter $\varepsilon\in \R_{>0}$,
which has $C^4(\bar\Omega\times \bar\Omega)$ classical regularity, its RKHS is a subset
of $C^2(\bar\Omega)$ and is norm-equivalent to
the Sobolev space $H^{(d+5)/2}(\bar\Omega)$.
However, in the following we will mostly focus on the infinitely smooth case,
especially use the Gaussian kernel with an isotropic
Euclidean distance $k^G_{\varepsilon}(x,x'):=\exp{(-\varepsilon^2\mynorm{x-x'}^2)}$,
or
involving an anisotropic Mahalanobis distance
\begin{equation}
  k^G_B(x,x'):=\exp{(-(x-x')^T B (x-x'))} \label{eqn:anisotropic-gaussian}
\end{equation}
for a symmetric positive definite matrix
$B\in \R^{d\times d}$.

\subsection{Generalized Interpolation}

We assume that we have chosen a kernel $k$ that matches the regularity of the solution,
i.e.\ $u\in \nsbar$. Then, we interpret the BVP as a collection of functional conditions:
We define sets of functionals via
$$
\Lambda_L :=\{\delta_x \circ L \, | \, x \in \Omega \},\quad
\Lambda_B:=\{\delta_x \circ B \, | \, x \in \partial \Omega \},
\quad \Lambda:= \Lambda_L \cup \Lambda_B,
$$
and corresponding target
values $y_\lambda:=\lambda(u) \in \R$ for $\lambda \in\Lambda$.
This means $y_\lambda=f(x)$ for
$\lambda= \delta_x \circ L \in \Lambda_L $
and $y_\lambda=g(x)$ if $\lambda= \delta_x \circ B \in \Lambda_B$.
Then the BVP can equivalently and compactly be expressed as
finding $u \in \nsbar$ such that
$$
   \lambda(u) = y_\lambda, \quad \lambda \in \Lambda.
$$
The set $\Lambda$ is linearly independent under the
uniform ellipticity condition.

The approximation now is obtained in the setting of generalized interpolation
\cite[Chap.\ 16]{wendland2005scattered}.
For this we assume to have chosen a set of functionals
$\Lambda_n:=\{\lambda_1,\ldots,\lambda_n\} \subset \Lambda$.
The corresponding trial space is set as
$$
  V_n := \myspan\{v_{\lambda_1},\ldots, v_{\lambda_n}\} \subset \nsbar,
$$
and the approximant $s_n \in V_n$ is obtained by requiring
$$
    \lambda_i (s_n) = y_{\lambda_i},\quad i=1,\ldots,n.
$$
Computationally, this corresponds to solving the following linear equation system
for $\alpha= (\alpha_i)_{i=1}^n \in \R^n$:
$$
K \alpha = y, \mbox{ with } K = \left(\scalprod{v_{\lambda_i}}{v_{\lambda_j}}_{\nsbar}\right)_{i,j=1}^n \in \R^{n\times n}
\mbox{ and } y=(y_{\lambda_1},\ldots,y_{\lambda_n})^T \in \R^n,
$$
and then setting
$s_n = \sum_{i=1}^n \alpha_i v_{\lambda_i} \in \nsbar$.
Since the set $\Lambda_n\subset\Lambda$ is linearly independent, this matrix $K$ is positive definite
thanks to the strict positive definiteness of the kernel. Thus a unique solution exists.
Moreover this approximant $s_n$ can be shown to be the solution of the
constrained norm-minimization problem
$$
   \min_{s \in \nsbar} \mynorm{s}^2_{\nsbar}\quad \mbox{ s.t. } \quad \lambda_i(s) = y_{\lambda_i},\ i=1,\ldots,n.
$$
Finally, the approximant $s_n$ can equivalently be characterized by an orthogonal
projection, i.e.\ 
$s_n= \Pi_{V_n}(u)$, where $\Pi_{V_n}:\nsbar \rightarrow V_n$ indicates the orthogonal
projection operator onto $V_n$.

For continuous linear functionals $\lambda \in \nsbar'$
one can show that the Riesz representer is given as the
functional applied to one argument of the kernel, i.e.
$v_\lambda = \lambda^{[1]} k(\cdot,\cdot)$, where the superscript 
indicates
application to the first argument such that the resulting object is a function only of the
second argument.
For the functionals of the BVP this results in
$v_\lambda= L^{[1]}k(x,\cdot)$ if $\lambda= \delta_x \circ L \in \Lambda_L$,
and similarly for $B$.

Assuming w.l.o.g.\ that the functionals in $\Lambda_n$ are sorted
into $n_L$ functionals from $\Lambda_L$ and $n_B$ from $\Lambda_B$
and the evaluation points are denoted as $x_{L,i},i=1,\ldots,n_L$ and
$x_{B,i},i=1,\ldots,n_B$,
the
kernel matrix can explicitly be written in block-matrix form as
$$
K = \left( \begin{array}{c|c}
             \left(L^{[1]}L^{[2]} k(x_{L,i},x_{L,j}) \right)_{i,j=1}^{n_L} &
                                                                             \left(L^{[1]}B^{[2]} k(x_{L,i},x_{B,j}) \right)_{i,j=1}^{n_L,n_B} \\
             \hline
          (*)^T   &  \left(B^{[1]}B^{[2]} k(x_{B,i},x_{B,j}) \right)_{i,j=1}^{n_B}                         
           \end{array}  \right),
$$
where the $(*)^T$ indicates the corresponding transpose
of the right upper block.
Due to the inherent symmetry of the system matrix, 
this PDE approximation procedure
is referred to as symmetric collocation. 

Moreover, the procedure recovers the standard kernel interpolation procedure for function approximation
in the case of assuming $L=\mathrm{Id}, B = \mathrm{Id}$, and hence
$\lambda_i = \delta_{x_i}$ for suitable
centers $x_i\in \bar \Omega$.

An important ingredient for the subsequent scheme is the generalized
power function $P_{\Lambda_n}:\nsbar' \rightarrow \R$ defined as
$$
P_{\Lambda_n}(\lambda) := \sup_{0\not = u \in \nsbar}
\frac{|\lambda(u-\Pi_{V_n}(u))|}{\mynorm{u}_{\nsbar}},\quad
\lambda\in \nsbar'.
$$
Observe that $P_{\Lambda_n}(\lambda)\leq \left\|\lambda\right\|_{\nsbar'}$ since $\Pi_{V_n}$ is an orthogonal projection, thus $P_{\Lambda_n}$ is uniformly bounded over $\Lambda$ if $\Lambda$ itself is bounded in $\nsbar'$.

The relevance of this quantity is twofold. First,
introducing $e_n:= u-s_n = u-\Pi_{V_n}(u)$ as the approximation error
of the collocation scheme, the definition of the power function
trivially yields an error bound for functional evaluations
$$
|\lambda(e_n)|=
|\lambda(u-\Pi_{V_n}(u))| \leq P_{\Lambda_n}(\lambda) \mynorm{u}_{\nsbar}, \quad \lambda \in \Lambda. 
$$
Assuming boundedness of $\Lambda$, the
well-posedness estimate
\eqref{eqn:well-posedness-estimate}
directly results in an $L^\infty$-error bound
\begin{eqnarray}
\mynorm{u-s_n}_{L^\infty(\bar \Omega)}
 &\leq & C_L \sup_{\lambda \in \Lambda_L} |\lambda(u-s_n)| 
  + C_B \sup_{\lambda \in \Lambda_B} |\lambda(u-s_n)|   \label{eqn:linfty-bound}
   \\
    &
   \leq & C'' \sup_{\lambda \in \Lambda} P_{\Lambda_n}(\lambda) \mynorm{u}_{\nsbar},
      \nonumber
\end{eqnarray} 
with $C'':=\max(C_L,C_B)$. 
Second, the power function is a measure of stability of the
(generalized) interpolation procedure, as will be explained
in the next subsection.

\subsection{Greedy Generalized Kernel Interpolation}

The number $n$ that can practically be considered in the full collocation approach
of the previous section is limited to at most a few thousand, due to several reasons.
From standard kernel interpolation it is well known that the
kernel matrix is in general dense, posing limitations due to
quadratic memory scaling and cubic complexity for the linear system solve.
Furthermore the matrix typically is ill-conditioned, which is getting more
severe with increasing $n$.
The same problems occur here, as the
system matrix is obtained by evaluations of
derivatives of the initial kernel, as exemplified in the previous section.
Although reducing the smoothness of the kernel generally improves the conditioning of its kernel matrix,
the kernel derivatives result in the scale parameters, e.g. $\varepsilon$,
occurring as factors, which may additionally lead to worse conditioning.

A method for counteracting these problems is to turn to
greedy procedures for selecting a subset of functionals.
For this we assume
to have a large set $\Lambda_N \subset \Lambda$ with $N\in \N$, for
which we cannot or do not want to
set up the full collocation matrix, hence we rather consider this to be
a set of candidate collocation functionals. We aim to incrementally select
functionals $\lambda_i \in \Lambda_N$ for $i=1,\ldots, n$,
with $n\ll N$, that provide good overall approximation by collocation on
a corresponding small set $ \Lambda_n \subset \Lambda_N$.

The abstract algorithm is given in Alg.\ \ref{alg:PDE-greedy}, which
requires a kernel $k$, the set of functionals $\Lambda_N$, the vector of
target values $y_N:=\Lambda_N(u) \in \R^N$ which means
the right-hand side (RHS) functions $f$ and $g$ evaluated in the corresponding collocation points,
as stopping criteria a maximum expansion size $n_{max}$, an accuracy tolerance
$\varepsilon_{acc}$ and a stability tolerance $\varepsilon_{stab}$.
The algorithm requires a general selection criterion
at iteration $i$
$$
  \eta_i:\Lambda_N \rightarrow \R_{\geq 0},
$$
which should indicate which functional is to be chosen next.
The criterion can make use of the current approximation space
$V_{i-1}$ and approximant $s_{i-1}$.
The specific choice of this selection criterion will be
given later.
The procedure then incrementally selects the
functional that maximizes the current selection indicator.
The corresponding Riesz-representer is computed,
orthonormalized with respect to the current space such that the next orthonormal basis vector
is obtained. The functional set and
trial space are correspondingly extended and the
next solution approximation computed by solving a (full) generalized kernel interpolation
problem {\tt{GeneralizedKernelInterpol}}
on the current subset $\Lambda_i$ as described in the
previous subsection. The notation $y_N(\Lambda_i)$ in this step
is to be understood as suitable subvector extraction.
These extension steps are iterated until any of the
stopping criteria is satisfied.

\begin{algorithm}
    \caption{\tt{GreedyGeneralizedKernelInterpol}}
    \label{alg:PDE-greedy}
    \begin{algorithmic}[1]
      \Require Kernel function $k$, set of functionals $\Lambda_N$,
      target values $y_N = \Lambda_N(u)$, maximum iterations $n_{max}$,
      accuracy tolerance $\varepsilon_{acc}$, stability tolerance $\varepsilon_{stab}$, selection criterion $\eta_i$ 
      \State $V_0 \gets \emptyset$, $\Lambda_0 \gets \emptyset$, $s_0 \gets 0$, $i \gets 0$,
         $P_0 \gets \infty$
        \While{$i \leq n_{max}$ \textbf{and} $\max |y_N - \Lambda_N(s_i)| > \varepsilon_{acc}$
             \textbf{and} $ P_{i}>\varepsilon_{stab}$ }
            \State $i \gets i+1$
            \State $\lambda_i \gets \arg\max\limits_{\lambda \in \Lambda_N} \eta_i(\lambda)$
            \State $P_i \gets P_{\Lambda_{i-1}}(\lambda_i)$
            \State $v_i \gets \frac{v_{\lambda_i} - \Pi_{V_{i-1}} v_{\lambda_i}}{P_i}$
            \State $\Lambda_i \gets \Lambda_{i-1} \cup \{\lambda_i\}$
            \State $V_i \gets V_{i-1} + \operatorname{span}(v_i)$
            \State $s_i \gets \texttt{GeneralizedKernelInterpol}(k, \Lambda_i, y_N(\Lambda_i))$
            \EndWhile
            \State $n \gets i$, 
            $s_n \gets s_i$, $\Lambda_n \gets \Lambda_i$\\
            \Return
            Iteration number $n$, approximant $s_n$, selected functionals $\Lambda_n$       
            \end{algorithmic}
\end{algorithm}

In the algorithm, the use of the power function
as a normalizing factor can be observed. If the power function
value is too small, this normalization can become too inaccurate.
Therefore the power function values are used as
an additional stopping criterion for preventing instability.

The overall selection criterion that we apply is
\begin{equation}
 \eta_i(\lambda) := w(\lambda) \, |\lambda(e_{i-1})|^\beta \, P_{\Lambda_{i-1}}(\lambda)^{1-\beta}, \label{eqn:beta-greedy-criterion}
\end{equation}
which depends on the current approximation error $e_{i-1}$, the
current generalized power function $P_{\Lambda_{i-1}}$,
a parameter $\beta \in [0,\infty]$ and a factor $w(\lambda)$.
We motivate the latter two ingredients: 
As the operators $L$ and $B$ in general do not coincide in their
differential order,
different powers of the shape parameter $\varepsilon$ appear in the corresponding matrix sub-blocks,
which especially result in potentially large variations in the magnitude of the
diagonal values.
Therefore, we introduce additional weighting factors $w(\lambda)$
for each functional
by which their importance can be adjusted. For instance, it is typically important to
ensure that sufficiently many Dirichlet functionals are chosen to
ensure the correct values in those points and not only select interior
differential functionals.

Further, as has been done in standard greedy kernel
interpolation \cite{wenzel2022analysis},
a parameter $\beta \in [0,\infty]$ is included in the
selection scheme, which allows to balance the importance of
the RHS data functions. Specifically the value $\beta=0$ is agnostic
of the RHS functions as it only is based on the power function. Hereby
it provides functionals that are good for all possible right hand sides,
but certainly provides suboptimal selection if only a single BVP with a
specific set of RHS data functions needs to be solved. In contrast, for
``target-dependency'', adaptivity to the RHS data functions can be turned
on by increasing $\beta>0$, especially covering the analogues of
the $f$-greedy ($\beta=1$), $f\cdot P$ ($\beta=1/2$) and
$f/P$-greedy ($\beta= \infty$) procedures of function interpolation. 
Using this indicator, we conclude from the algorithm that, as long as the loop is continued -- i.e., the accuracy threshold is not satisfied -- there exists a functional for which an interpolation error occurs,
hence $\eta_i(\lambda) > 0$ for at least one functional.
A reasonable property of the selection criterion now is that
$\eta_i(\lambda_j)=0$ for all $j=1,\ldots,i-1$ as the
(generalized) interpolation conditions are met at the selected functionals.
This guarantees that no previously selected functional is chosen twice.

As the $\beta$-greedy indicator and the iterative approximation scheme
do not require the full kernel matrix, the overall scheme is very efficient.
In linear algebra terms it corresponds to a partial pivoted Cholesky
decomposition of the system matrix without assembling it in full,
but only iteratively computing single rows. The power function
values
$P_i$ of step 5 in Alg.\ \ref{alg:PDE-greedy} are the diagonal elements of the Cholesky factor, further explaining their role
in controlling the stability of the algorithm.
The pivoting rule is not one of the typical criteria used in
numerical linear algebra but exactly corresponds to the
$\beta$-greedy criterion, except for $\beta=0$, which corresponds to a classical full pivoting.
The Gram-Schmidt orthonormalization that is performed inside the algorithm
generates an orthonormal basis of $V_n$ that is the generalization
of the so called Newton-basis in kernel interpolation.
Especially the incremental approximants can be efficiently updated, as
the previous coefficient vector can be reused and only
a single new coefficient for the new basis vector needs to be computed
for the generalized interpolation. 

Overall, the training complexity is then $\cO(Nn^2)$, i.e.\
only linear in $N$, hence allowing $N$ in the order of millions on a
standard PC.
The prediction complexity at one point is only linear $\cO(n)$, hence very fast
when assuming a small expansion size $n$. 

The overall greedy generalized kernel interpolation method using the selection indicator
\eqref{eqn:beta-greedy-criterion} is denoted as generalized $\beta$-greedy method.
If $\Lambda$ is given by point evaluations from
PDE BVPs, we denote the approach as PDE-$\beta$-greedy method.

\begin{rem}[Interpretation as Abstract Plain Kernel Interpolation]
  If the generalized kernel interpolation procedure is well-posed, especially
  uniqueness for given RHS functions holds, then the generalized kernel interpolation
  procedure can be shown to be equivalent to a standard
  plain kernel interpolation problem on the set $\Lambda$ using the
  functional kernel
  $k_{\Lambda}(\lambda,\lambda'):= \scalprod{v_{\lambda}}{v_{\lambda'}}$.
  Moreover, since an equivalence of the corresponding power functions can also be shown,
  the corresponding generalized $\beta$-greedy method corresponds to an abstract
  plain $\beta$-greedy interpolation procedure.
  For details, we refer to \cite{HS25}.
  The relevance of this correspondence is that plain $\beta$-greedy kernel interpolation code can
  be reused for approximating BVP by just replacing the involved kernel.
\end{rem}

\begin{rem}[Relation to GEIM Method]
  The generalized empirical interpolation method (GEIM)
  \cite{Mula2016} on first sight may seem similar to greedy
  generalized interpolation,
  as it also is based on selecting functionals for providing approximants.
  In general, however, the setting is rather different: 
  The GEIM aims at approximating any function from a given set $F$ of functions rather than approximating a single target function.
  Also
  the selected basis functions in GEIM are functions from the target set instead of
  Riesz-representers of functionals.
  However, there is a special artificial situation where the methods become comparable:
  If the target set is $F=\{k(x,\cdot)\ | x\in \Omega\}$ for a positive
  definite RBF kernel and the functionals are merely point evaluations,
  then the functionals selected by the generalized $P$-greedy method have Riesz-representers
  which are elements from $F$ and exactly correspond to the basis chosen by the GEIM method.
  This equivalence transfers to the PDE-$P$-greedy procedure, if the abstract
  kernel $k_\Lambda(\lambda,\lambda'):=\scalprod{v_\lambda}{v_{\lambda'}}$
  can be written as radial kernel depending on the distance of the sampling points
  corresponding to $\lambda$ and $\lambda'$. Indeed, the fact that $k_\Lambda$ is radial ensures that it is uniquely maximized at 
$\lambda=\lambda'$.  
  For other values than $\beta=0$ there cannot be an equivalence to GEIM, as the
  (PDE-)$\beta$-greedy methods for $\beta>0$ provide adaptivity to each single given target function.  
\end{rem}

\subsection{Convergence Rates}
\label{sec:convergence-rates}
For the case of finite smoothness of the target function, i.e.\ assuming a
kernel $k$ for which the RKHS $\nsbar$ is norm-equivalent to some
Sobolev-space, algebraic convergence rate results can be obtained.
First, for the case of second order elliptic BVP functional evaluations 
we state a simplified (but coarser) bound that can be obtained
from \cite[Thm 5.1]{PDEgreedy2025}.
Especially we modify the index range and upper bound a
residual term as presented in the proof in the appendix.

\begin{prop}[Algebraic Convergence for BVP Residuals]
\label{prop:algebraic-rates-functional-residuals}
Assume $k$ is a kernel on $\bar \Omega\subset \R^d$
such that $\nsbar \asymp H^{\tau}(\bar \Omega), \tau>d/2 +2$.
If $u\in \nsbar$, then any
PDE-$\beta$-greedy scheme for $\beta\in [0,1]$ with uniform weighting $w(\lambda)=1$
provides an approximation
sequence $s_n$ such that for $n \geq 6$ and
abbreviating $ \alpha:=(\tau-2)/d-1/2 >0$ it holds
  \begin{equation}
  \min_{i=1,\ldots,n} \sup_{\lambda \in \Lambda} |\lambda(u-s_i)|
  \leq C n^{-\beta/2} n^{-\alpha} \log(n)^\alpha \mynorm{u}_{\nsbar}.
  \label{eq:aux}
  \end{equation}
\end{prop}
Analogously to Thm.~5.1 in \cite{PDEgreedy2025},
applying the well-posedness bound for the BVP problem
\eqref{eqn:well-posedness-estimate}
directly yields the same rate for the error in the $L^\infty$ norm with a modified constant.

\begin{cor}[Algebraic Convergence Rate for the $L^\infty$-Error]
    \label{cor:algebraic-rates-linfty-error}
  Under the assumptions (and especially the constant $C$)
  of Prop.\ 
\ref{prop:algebraic-rates-functional-residuals}
we have
  \begin{equation}
  \min_{i=1,\ldots,n} \mynorm{u-s_i}_{L^\infty(\bar \Omega)}
  \leq C \max(C_L,C_B) n^{-\beta/2} n^{-\alpha} \log(n)^\alpha \mynorm{u}_{\nsbar}.
  \end{equation}
\end{cor}
Note that similar statements can be formulated
for $\beta\in (1,\infty]$, cf.\ \cite{PDEgreedy2025}.
The striking feature of such bounds is that the $n^{-\beta/2}$ factor
is independent of the dimension and especially indicating that
$\beta>0$, i.e.\ schemes that involve the target values for
functional selection have a better convergence rate
than the case $\beta=0$, which is agnostic of the RHS
functions $f,g$.

For the case of infinitely smooth kernel and solution, exponential
convergence can be proven, cf.\ Thm.\ 3.18 in 
\cite{vogel24}, which we slightly simplify similarly to
Prop.\ \ref{prop:algebraic-rates-functional-residuals}, as presented in
the proof in the appendix. 

\begin{prop}[Exponential Convergence for BVP Residual]
\label{prop:exponential-rates-functional-residuals}
Let $k$ be the Gaussian kernel on a polygonal $\bar \Omega\subset \R^d$ and
$\beta=1$.
If $u\in \nsbar$, then the PDE-$\beta$-greedy scheme with uniform weighting
$w(\lambda)=1$ provides an approximation
sequence $s_n$ such that for $n\geq 6$
  \begin{equation}
  \min_{i=1,\ldots,n} \sup_{\lambda \in \Lambda} |\lambda(u-s_i)|
  \leq C n^{-\beta/2} e^{-c_1 n^{\frac{1}{2d}}} \mynorm{u}_{\nsbar}.
  \label{eq:aux3}
\end{equation}
\end{prop}

Again, using the BVP well-posedness estimate
\eqref{eqn:well-posedness-estimate}
we conclude the same
exponential rate but with modified constants for the $L^\infty$-error, which is a simplification
of Cor.\ 3.19 in \cite{vogel24}.

\begin{cor}[Exponential Convergence for $L^\infty$-Error]
  \label{cor:exponential-rates-linfty-error}
  Under the assumptions (and especially the constant $C$)
  of Prop.\ 
\ref{prop:exponential-rates-functional-residuals}
we have
  \begin{equation}
  \min_{i=1,\ldots,n} \mynorm{u-s_i}_{L^\infty(\bar \Omega)}
  \leq C \max(C_L,C_B) n^{-\beta/2} e^{-c_1 n^{\frac{1}{2d}}} \mynorm{u}_{\nsbar}.
  \end{equation}
\end{cor}

We want to emphasize here a clear and striking advantage of these
results over neural networks:
These convergence rates are not only
expression rates (i.e.\ existence statement for
corresponding approximant with specified accuracy), but
as well constructive, i.e.\ after the stated number of iterations
a corresponding approximant is guaranteed to be computed by the greedy strategy.
  For neural networks, however, in general no guarantees can be given that the
  optimization procedure actually finds a suitable approximant in a limited number
  of iterations.

We further want to relate to convergence statements on Gaussian-Process
based PDE-approximation \cite{OwhadiErrorAnalysis} which 
states algebraic convergence in the so called fill-distance $h_X$ of
the points used for collocation. Hence there is only a convergence
guarantee if the collocation points get dense in $\Omega$.
In contrast, the $\beta$-greedy approximation can provide
convergence statements only in terms of the expansion size $n$ and
avoiding the fill-distance, hence these results
do not require that the points get dense everywhere in $\Omega$.
This is considered to be an conceptional advantage of the greedy approach. 
   
\section{Extension to the Parametric Case}

\label{sec:muPDEgreedy}

We now formulate the parametric extension of the problem to be considered and the approximation technique.
We assume to have a bounded parameter set $\cP \subset \R^{d_\mu}$
and a possibly parametrized geometry, i.e.\
for all $\mu\in\cP$ we have an (open) domain $\Omega(\mu) \subset \R^{d_x}$ for the
position variable $x$, and indicate its closure as $\overline \Omega(\mu)$. We
set the dimension of the overall position-parameter domain as $d:=d_x + d_\mu$. 
We define the set $\bOmega := \{ \bx= (x,\mu) \in \R^d \, | \, \mu \in \cP, x\in \Omega(\mu) \}$
for position-parameter tuples $\bx=(x,\mu)$.
As the differential operator only acts on the boundary points of $\Omega(\mu)$ and not on
the overall boundary of $\bOmega$, we define
$\partial \bOmega := \{ (x,\mu) \in \R^d \, | \, x\in \partial \Omega(\mu),\mu \in \cP \}$.
In the simplest case of $\Omega(\mu)= \Omega$ being constant w.r.t.\ $\mu$, this boundary set
corresponds to a cylinder mantle in the position-parameter domain. 
Subsequently we define the union of these sets 
as
$\overline \bOmega := \bOmega \cup \partial \bOmega$. Observe that this is just a notation in analogy to the non-parametric case, and the sets are not necessarily the interior, boundary, and closure of a given set $\bOmega$
(e.g.\ in the given example $\partial \bOmega$ does not contain the cylinder top/bottom faces).
In particular, this general formulation allows discrete parameter
sets $\cP$ without assuming any common structure between the PDEs
corresponding to different values of $\mu$.

We assume that for any $\mu\in \cP$  a corresponding parametric PDE BVP problem is given, i.e.\
find $u(\cdot,\mu)\in C^2(\Omega(\mu))\cap C^0(\overline \Omega(\mu))$ such that
$$
L(\mu) (u(\cdot,\mu)) = f(\mu) \mbox{ in } \Omega(\mu), \mbox{ and }
B(\mu) (u(\cdot,\mu)) = g(\mu) \mbox{ on } \partial \Omega(\mu).
$$
Here $L(\mu)$ is assumed to be a second order operator of the form
$$
(L(\mu) u)(x) = \sum_{i,j=1}^{d_x} a_{ij}(x,\mu) \partial_{x_i} \partial_{x_j} u(x)
   + \sum_{i=1}^{d_x} b_i(x,\mu) \partial_{x_i} u(x) + c(x,\mu) u(x),
$$
with coefficient functions $a_{ij},b_i,c: \bOmega \rightarrow \R$,
and similarly for the boundary value operator $B(\mu)$.
We assume that the problem is well-posed, i.e.\ for any $\mu$ we assume the existence of a unique classical solution $u(\cdot,\mu)$.
Thus, the parametric solution can be understood as a
function $u:\overline \bOmega \rightarrow \R$.
Furthermore we suppose a certain regularity of the solution
in the position and the parameter, i.e.\ we assume that
$u \in \calh_k(\overline \bOmega)$ for some position-parameter kernel
$k:\overline \bOmega \times \overline \bOmega \rightarrow \R$.
For simplicity we restrict to
product kernels, i.e.\ $k(\bx,\bx') = k_x(x,x') \cdot k_\mu(\mu,\mu')$
for suitable position kernel $k_x:\R^{d_x} \times \R^{d_x} \rightarrow \R $ and parameter kernel
$k_\mu: \R^{d_\mu} \times \R^{d_\mu} \rightarrow \R$.
Such product kernels correspond to tensor products of RKHSs

 and are very versatile \cite{Aronszajn1950,Albrecht2025}. 
Then $\calh_k(\overline \bOmega)$ is a restriction of the tensor product space
$\calh_{k_x}(\R^{d_x}) \otimes \calh_{k_\mu}(\R^{d_\mu})$ to the domain $\overline \bOmega$.
The advantage of this structure is that the differential operator
only acts on the position variable, for instance
$L(\mu)^{[1]}k(\bx,\bx') = (L(\mu)^{[1]} k_x (x,x')) \cdot k_\mu(\mu,\mu')$. 
This simplifies the practical realization as generic and modular implementation and combination
of the corresponding kernel functions is possible. E.g., different regularity assumptions w.r.t.\
the parameter and the
position could be realized by choosing corresponding position and parameter kernels. 
Assuming that $f,g$ and the coefficients $a_{ij},b_i,c$ are $C^\infty$  with respect
to the position variable and $\partial \Omega(\mu)$ is a $C^\infty$-boundary,
uniform ellipticity of $L(\mu)$ and $c\leq 0 $ implies existence, uniqueness and
that $u(\cdot,\mu)$ is infinitely smooth with respect to the position variable. 
From the RB literature it is known that differentiability of coefficients w.r.t.\  the parameter
transfers to corresponding differentiability of the solution with respect to the parameter
\cite{Haasdonk2017}.
But such smoothness in the solutions can occasionally as well be obtained despite
less regular data or geometry. As an example we will consider a non-Lipschitz
domain including topology change in the experiments.

The generalized interpolation setting can then directly be applied by
identifying
$$\Lambda:=\Lambda_L\cup \Lambda_B:=\{\delta_x\circ L(\mu) \, | \, (x,\mu)\in \bOmega\}
\cup \{\delta_x \circ B(\mu) \, | \, (x,\mu) \in \partial \bOmega\} \subset
\calh_k(\overline \bOmega)',$$
and applying the method from Sec.\ \ref{sec:PDEgreedy}
using the position-parameter domain $\overline \bOmega$ instead of the
position domain $\bar \Omega$. Especially the final approximant
is $s_n:\overline \bOmega \rightarrow \R$, which can be interpreted as
a global a-priori surrogate model similar to PINNs \cite{PINNref}, PGD \cite{nouy2010}, etc.
The ``offline-phase'' consists of a single training phase of the model
(sampling of the position-parameter domain and boundary, greedy search of collocation
points and computation of approximant). The ``online-phase'' then simply consists of
pointwise evaluation of the surrogate in any new position-parameter
tuple $\bx\in \overline \bOmega$.

We want to explicitly state that we avoid a ``standard approach'' of parametric
approximation, which would consist of providing some solutions (snapshots) for some
fixed parameter samples. Then an overall
parametric approximant would be obtained by some sort of parametric interpolation or
other linear combination of these precomputed particular solutions. 
The reason for not following this approach is the overall solution expansion size, which would
grow multiplicatively in the number of parameter samples used for the parametric interpolation
(or more general linear combination) and the individual expansion size of each solution snapshot.
We expect that by avoiding this tensor-product construction, greedy strategies over the
parameter-state domain are able to place samples more intelligently and sparsely
than choosing a fixed tensor-product grid.

We now want to briefly discuss the changes in the error analysis for the parametric extension.
First, with respect to the position-parameter domain, the elliptic (w.r.t. position variable)
differential operator
is degenerate in the parametric direction as there are no parametric derivatives.
So the framework of \cite{PDEgreedy2025} might seem not fully suitable at first glance.
However, when inspecting the proof of the convergence rates leading to
Prop.\ \ref{prop:algebraic-rates-functional-residuals} or
Prop.\ \ref{prop:exponential-rates-functional-residuals}, one realizes that actually
no ellipticity assumption is required for the convergence rates of the functional
residuals. Only for the conclusion of the $L^\infty$-error convergence rates
of Cor.\
\ref{cor:algebraic-rates-linfty-error} and Cor.\ \ref{cor:exponential-rates-linfty-error},
a well-posedness result is required. If we assume uniform (in parameter) ellipticity
of the differential operator and uniform (in parameter) boundedness of
its coefficients and uniform boundedness of the RHS data functions, Eqn.\ 
\eqref{eqn:well-posedness-estimate} can be formulated with constants $C_L,C_B$
which are independent of $\mu \in \cP$.
Then, especially, the bound \eqref{eqn:linfty-bound} 
holds for the whole position-parameter domain, providing an upper bound for the error on single
slices $\Omega(\mu)\times\{\mu\}$, hence for any $\mu\in \cP$ it holds
\begin{eqnarray}
  \mynorm{u(\cdot,\mu)-s_n(\cdot,\mu)}_{L^\infty(\bar \Omega(\mu))}
& \leq &
\mynorm{u-s_n}_{L^\infty(\bar \bOmega)} \nonumber \\
  &
  \leq & C_L \sup_{\lambda \in \Lambda_L}{\lambda(u-s_n)} +
  C_B \sup_{\lambda \in \Lambda_B}{\lambda(u-s_n)}. 
\label{eqn:parametric-well-posedness-estimate}
\end{eqnarray}
Then, also convergence rates for the RHS transfer to convergence
rates for the left hand side (LHS) restricted to slices.

\section{Numerical Experiments}

\label{sec:experiments}

The following experiments were performed with MATLAB R2024b
on a system with i7-1355U 1700 MHz CPU and 32 GB memory. 
For full reproducibility we
provide the code\footnote{\url{https://gitlab.mathematik.uni-stuttgart.de/pub/ians-anm/mu-pde-beta-greedy}.}.

In all examples, we consider a standard Poisson problem with Dirichlet
boundary values, i.e.\ $L = -\Delta,\ B= \mathrm{Id}$, hence we only
need to specify the domains $\cP, \Omega(\mu)$ and the RHS data
functions $f(\cdot,\mu),g(\cdot,\mu)$
in the following model problems. 

\subsection{Nonaffine Geometry Parametrization}

\label{sec:moving-circles}
We start with an example of a non-affine geometry parametrization of ``moving circles'',
where the position-domain is $\Omega(\mu) = B_1(0)\cup B_1((\mu,0)^T) \subset \R^2$ 
for $\mu\in[0,1.1]=:\cP \subset \R$, hence $d_x=2, d_\mu = 1$ and the overall state dimension is
$d=3$. The parameter is chosen such that the topology of the geometry changes from
a single connected to two disconnected circular components. The kinks in the domain boundary
for $\mu\in(0,1)$ are such that the domain $\Omega(\mu)$ is
non-convex. Such a geometry parametrization is clearly non-affine
and no fixed spatial point set can be chosen, and no standard RB or any 
parametric interpolation or linear combination of spatial ``snaphot'' solutions would be suitable
as such would not exactly recover the boundary geometry.

\begin{figure}[h!t]
\begin{center}
 a) \includegraphics[width=0.4\textwidth]{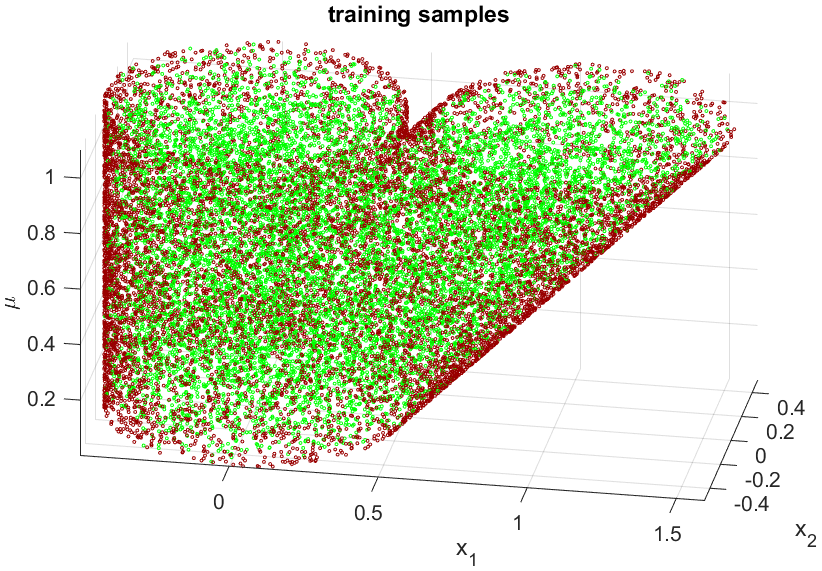} $\quad$ b)
\includegraphics[width=0.4\textwidth]{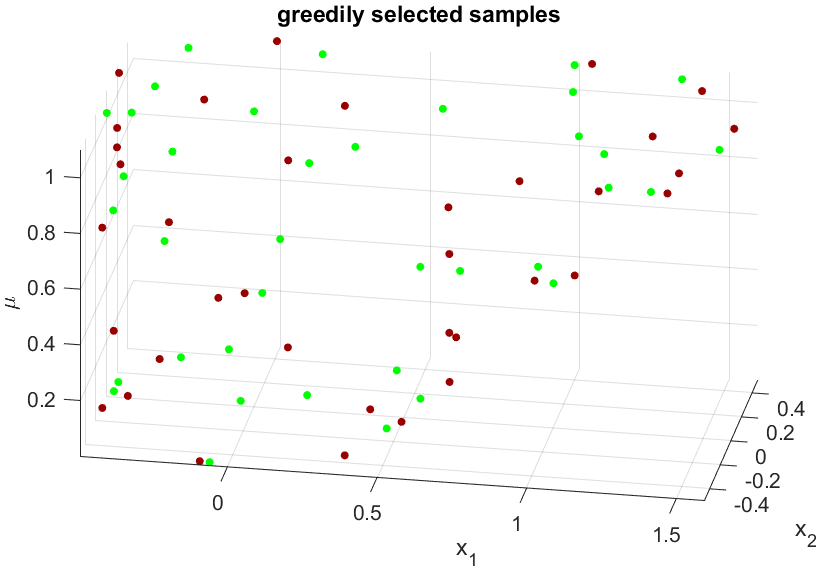}  
\end{center}
\caption{\label{fig:moving-circles-sampling}
  Functional samples for the moving circles example, indicated by
  point evaluation sites for the interior (green) and
  boundary operator (red). Subplot a) indicates the training samples. Subplot b) illustrates
  the greedily selected functionals for the smooth solution case for selected hyper-parameters.}
\end{figure}

The position-parameter domain $\bar \bOmega$ is consisting of the union of two (skew) cylinders and as
training set of the functionals we sample randomly and uniformly $N_L:=10^4$ interior
and $N_B:=10^4$ boundary points, hence
$N=N_L+N_B=2\cdot 10^4$ sampling points, which determine the corresponding
training set of functionals, by corresponding point evaluations of either
the differential or boundary operator, and are visualized in Fig.\ \ref{fig:moving-circles-sampling}a).

We consider two cases, first a case of a smooth solution, then a more challenging
setting with a singular solution.

For the smooth solution case, we assume the
known solution $u(x,\mu) = \frac12 (\mynorm{x}^2 + \mu^2)$
and choose $f:= L u$, $g:=u|_{\partial \bOmega}$ 
correspondingly.
This allows us to assess true pointwise errors as quality criterion.

Anticipating the $C^\infty$ regularity of the solution, we choose as kernel a product of
Gaussian kernels for both the
position and the parameter coordinate, but with different (squared) shape parameters
$\varepsilon_\mu^2,\varepsilon_x^2$. 

We choose an $f$-greedy strategy ($\beta=1$) as this seems most promising from the
convergence rate analysis.
We fix the maximum expansion size as
$n_{max}= 100$ and the stopping tolerances as $\varepsilon_{acc}= \varepsilon_{stab}= 1\cdot 10^{-15}$.
The functional weightings $w(\lambda)$ are chosen as $w_B \in \R_{>0}$ for $\lambda\in \Lambda_B$ (to be selected
later), and $w_L=1$ for $\lambda\in \Lambda_L$.

The (squared) shape parameters $\varepsilon_x^2,\varepsilon_\mu^2$ and the functional weighting parameter
$w_B$ are dominantly influencing the prediction accuracy of the models, hence are
considered to be hyper-parameters in this example. 
Therefore each of these parameters was consecutively 
selected by a 1D grid search over some discrete values of the parameter and
choosing the parameter resulting in the least ``loss'' of the approximant
measured over a randomly drawn validation set of 10000 interior $\Lambda_{L,val}$ and 10000 boundary
functionals $\Lambda_{B,val} $.
As loss function for validation of a model function we choose
a weighted sum of the maximum interior and maximum boundary BVP residual error
\begin{equation}
e_{val}(s_n) := \gamma_L \sup_{\lambda \in \Lambda_{L.val}} |\lambda(u-s_n)| +
             \gamma_B \sup_{\lambda \in \Lambda_{B,val}} |\lambda(u-s_n)|, \label{eqn:validation_loss}
\end{equation}
which is motivated by the RHS of the error equation \eqref{eqn:parametric-well-posedness-estimate}.
In lack of true values for the constants, we simply set $\gamma_L = \gamma_B = 1$.                          
The final selected parameters 
are $\varepsilon_x^2 = 0.0063096, \varepsilon_\mu^2 = 0.01, w_B=1.0$.

The greedily selected functionals (indicated by colors and their evaluation point location)
are given in Fig.\ \ref{fig:moving-circles-sampling}b), in which 38 interior and 35 boundary
functionals were chosen. 
One can clearly see the sparsity compared to the full training set of functionals in subplot a).
The greedy strategy seems to select the functionals in a well-distributed fashion.
The training time for selecting the 73 functionals and computing the approximate solution
was 0.32 seconds.

\begin{figure}[h!t]
\begin{center}
  a) \includegraphics[width=0.35\textwidth]{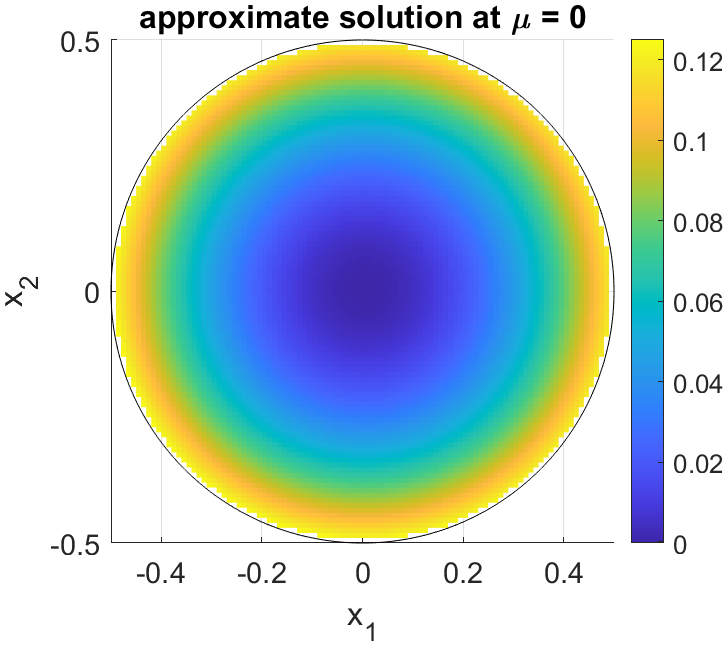}$\quad$
  b) \includegraphics[width=0.50\textwidth]{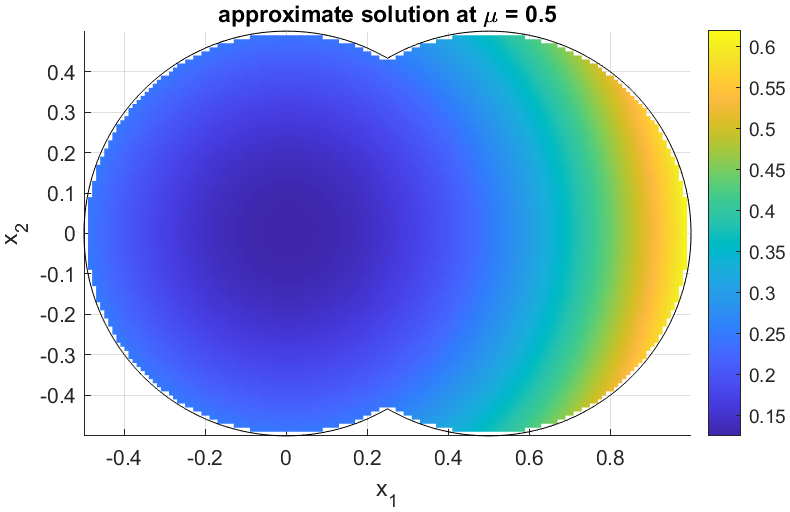} \\
  c) \includegraphics[width=0.7\textwidth]{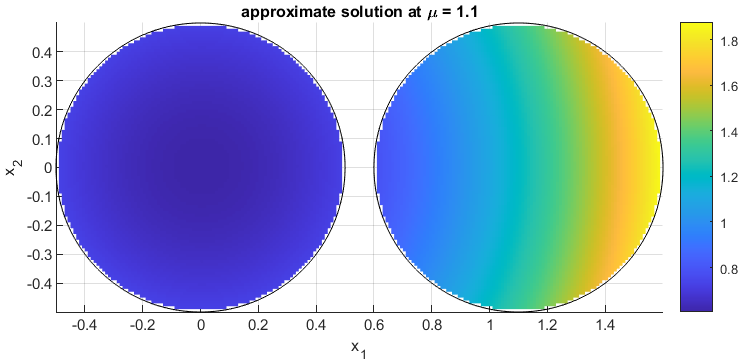}
\end{center}
\caption{\label{fig:moving-circles-smooth-solution}
  Visualization of the smooth solution case for 
  the moving circles example: Approximate solution for selected hyper-parameters at
  $\mu=0.0, 0.5$ and $1.1$.
}
\end{figure}

In Fig.\ \ref{fig:moving-circles-smooth-solution} we illustrate the parametric approximate
solution evaluated for several parameters. The changing topology is clearly visible.
The resulting absolute accuracy of the model is $\mynorm{u-s_n}_{L^\infty(\bar \bOmega)}=2.9719\cdot 10^{-7}$
measured over an independently drawn test set in $\bar \bOmega$.

\begin{figure}[h!t]
\begin{center}
  \includegraphics[width=0.6\textwidth]{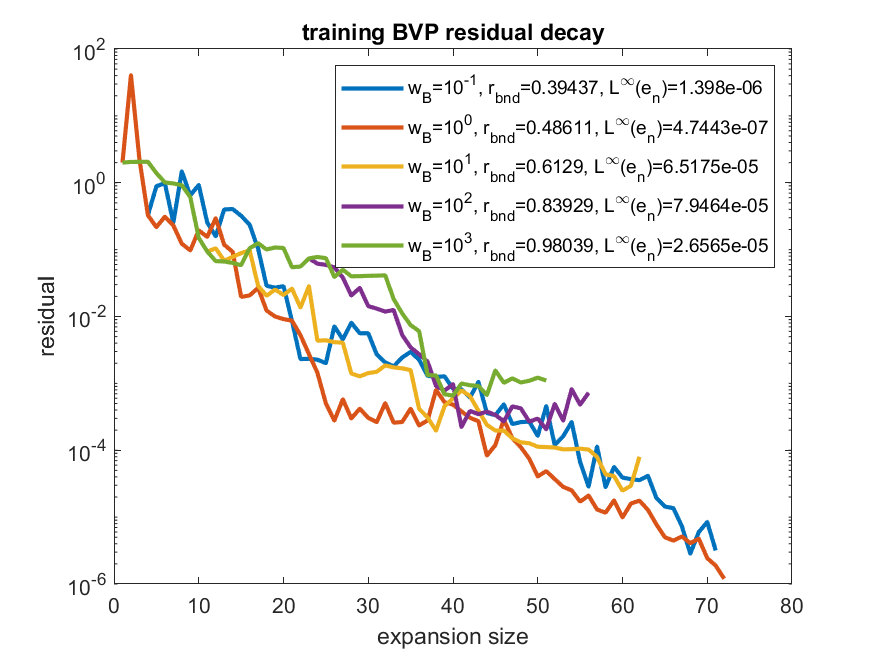}
\end{center}
\caption{\label{fig:moving-circles-convergence}
  Visualization of BVP training residual convergence, ratio of selected boundary functionals (indicated as $r_{bnd}$)
  and true $L^\infty$-error (indicated as $L^\infty(e_n)$) at maximum expansion size 
  for different boundary functional weighting values $w_B \in \{10^{-1},\ldots,10^3\}$
  for the moving circles example with the smooth solution case. 
}
\end{figure}

In Fig.\ \ref{fig:moving-circles-convergence} we illustrate the training BVP residual
decay, i.e. $\max_{\lambda \in \Lambda_N} |\lambda(u-s_n)|$.
One can clearly see the overall asymptotic exponential convergence,
which is what we expect in view of the error convergence analysis.
But also, we recognize local non-monotonicity or plateaus, which is why
the minimum in the LHS of the error bounds in Sec.\ \ref{sec:convergence-rates}
indeed is not an artifact but seems useful.
For each curve, we also report the ratio $r_{bnd} $ of selected boundary functionals over the total
number of selected functionals, which is
naturally increasing with growing $w_B$. Also, occasional early stopping of the greedy loop
due to numerical accuracy reasons can be observed. The $L^\infty$-errors of the final approximants
are as well reported
for each curve in the legend, 
and we note that indeed a local minimum is to be expected in the interval
$w_B \in [10^{-1},10^1]$, which is consistent with the optimized value $w_B=1$ reported earlier.

For the singular solution case, we choose $g(x,\mu)= 0$ and $f(x,\mu)=-1$.
Despite the simple data functions, the
geometrical parametetrization results in a lack of an
explicit closed form representation of the solution for the parameters $\mu \in(0,1)$.
For $\mu=0$ and $\mu\in[1,1.1]$ the solution $u(\cdot,\mu)$ is expected to be
$C^\infty(\Omega(\mu))$ due to the smoothness of the
data functions and the boundary $\partial \Omega(\mu)$.
But for $\mu\in(0,1)$ it is expected to have singularities in its
derivatives due to the entering corner. 
Due to this anticipated reduced regularity now the quadratic Mat\'ern
kernel \eqref{eqn:matern} is chosen.
Still, this clearly is an example where the target function $u$ is outside the RKHS $\calh_k(\bar \bOmega)$,
hence we do not expect the convergence rates as specified in Sec.\ \ref{sec:convergence-rates} and actually
the convergence is theoretically unclear at all. But still the example is interesting to assess.
We cannot access the true error of the approximants, but can only refer to the
BVP functional residual values.

Instead of $f$-greedy we turn to the $P$-greedy ($\beta=0$) strategy, as the target-data dependent criteria
tend to cluster points and imply early instability in the case of target functions outside the RKHS. 
In contrast,
the $P$-greedy strategy is known to produce more uniformly distributed points \cite{wenzel2021novel}. 
We fix the maximum expansion size as
$n_{max}= 1000$ and the stopping tolerances as before to $\varepsilon_{acc}= \varepsilon_{stab}= 1\cdot 10^{-15}$.
The shape parameters were manually chosen as $\varepsilon_x = \varepsilon_\mu = 3$ and
the functional weightings as $w_B=100000^{1/2}$ and $w_L=1$.

\begin{figure}[h!t]
  \begin{center}
    a) \includegraphics[width=0.5\textwidth]{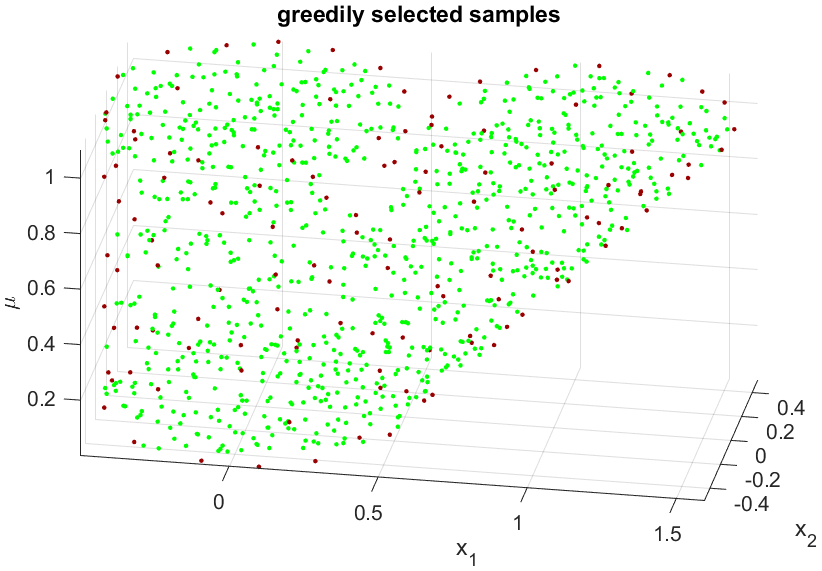}  \\
    b) \includegraphics[width=0.4\textwidth]{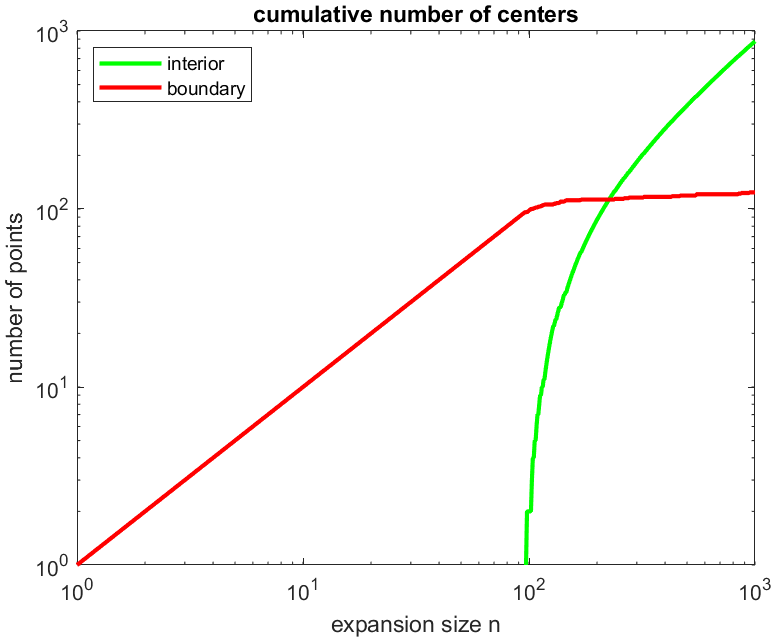} $\quad$  
    c) \includegraphics[width=0.4\textwidth]{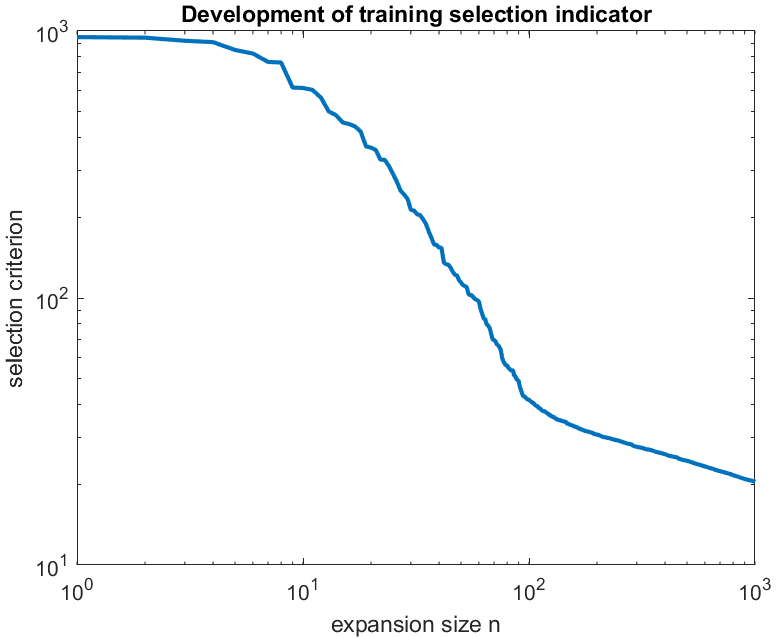}   
\end{center}
\caption{\label{fig:moving-circles-sampling-nonsmooth}
  a) Greedily selected functional samples for 
  the moving circles example with the
  nonsmooth solution case, indicated by
  point evaluation sites for the differential (green) and
  boundary operator (red).
  b) Cumulative number of boundary and interior functionals during greedy selection. 
  c) Decay of the greedy selection indicator during training. 
}
\end{figure}

In Fig.\ \ref{fig:moving-circles-sampling-nonsmooth}a)
we illustrate the
selected 876 interior and 124 boundary functionals. 
The uniformity and the considerably larger number of functionals compared to
the smooth example are evident. 
The training for selecting the 1000 functionals took 9.35 seconds.
This is approximately a factor 8 less than we would expect from a pure $n^2$ scaling compared to the
runtime in the smooth example. 
We report the cumulative number of boundary and interior functionals during the greedy selection in
Fig.\ \ref{fig:moving-circles-sampling-nonsmooth}b). 
The resulting maximum greedy selection criterion over the iterations is plotted in
Fig.\ \ref{fig:moving-circles-sampling-nonsmooth}c).
We nicely see a monotonic decay over
almost two orders of magnitude.
We can recognize two phases, where first only boundary functionals are chosen, then
dominantly the interior is being sampled.
This is reflected in the selection criterion (i.e. weighted power function) decay in different
decay regimes. The phase of predominantly boundary functional selection for $n$ up to 100 corresponds to
the pure (non PDE) function approximation of the infinitely smooth zero-boundary value function
in an effective $d-1$ dimensional
manifold with correspondingly expected more rapid decay.
When approximating the interior domain starting at about $n>100$,
the decay rate is reduced due to the higher degree of the differential operator, the full 
dimension $d_x$ of $\Omega(\mu)$ and the lower regularity of the target function.
The merely algebraic decay is very clearly visible.

\begin{figure}[h!t]
\begin{center}
  a) \includegraphics[width=0.35\textwidth]{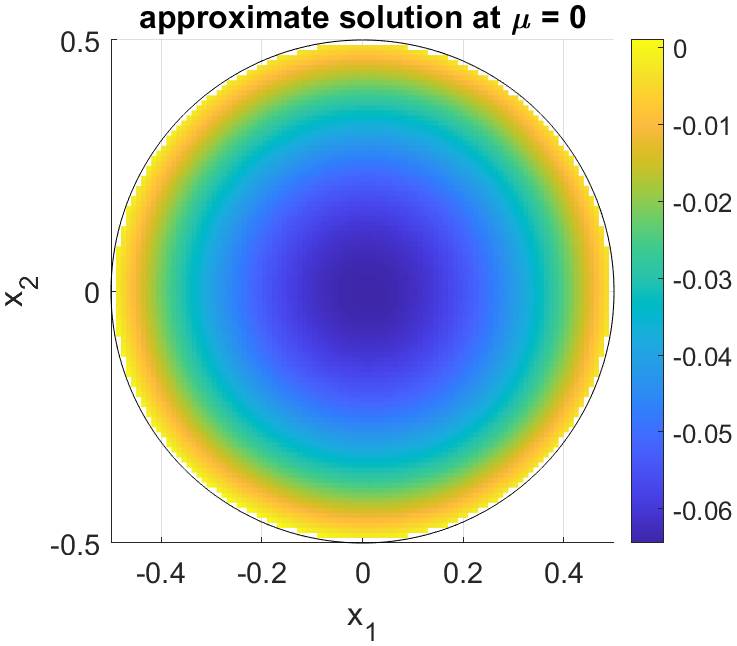}$\quad$
  b) \includegraphics[width=0.50\textwidth]{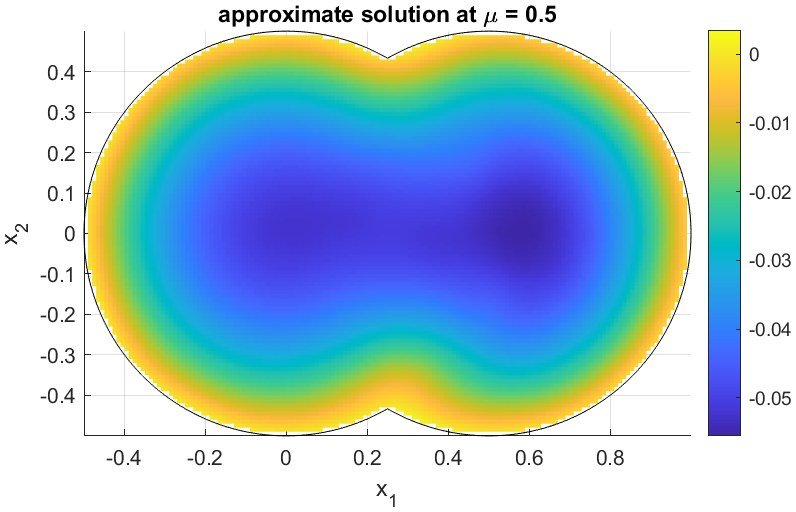} \\
  c) \includegraphics[width=0.7\textwidth]{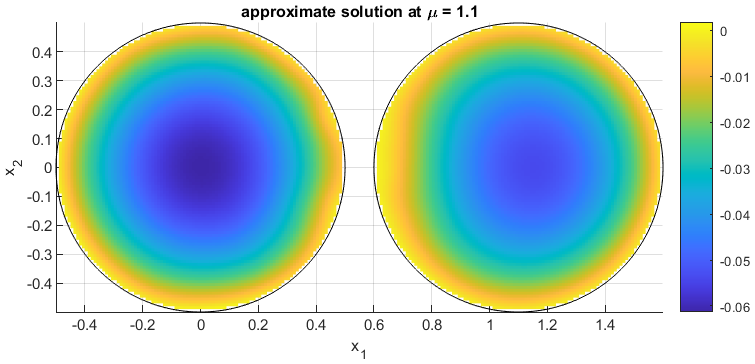}
\end{center}
\caption{\label{fig:moving-circles-nonsmooth-solution}
  Visualization of the 
  nonsmooth solution case for 
  the moving circles example, approximate solution at $\mu=0.0, 0.5$ and $1.1$.
}
\end{figure}

However, the difficulty of the given example is reflected by several observations. 
In Fig.\ \ref{fig:moving-circles-nonsmooth-solution} we illustrate the parametric approximate
solution evaluated for several parameters. As we have no closed form solution of the problem, we do not give
absolute error quantification. But qualitatively the solution seems reasonable in the sense that the boundary
values are visually accurately captured, and quantitatively realize a maximum absolute value over the boundary
training points of $ 0.00868$. Also the sign of the solution is consistent with the (negative) source
function. 
Due to the nonsymmetric point/functional placement, the pure expected rotational symmetry of the
solution for $\mu=0$ and $\mu=1.1$ is slightly distorted in the approximation. 

\begin{figure}[h!t]
\begin{center}
  a) \includegraphics[width=0.35\textwidth]{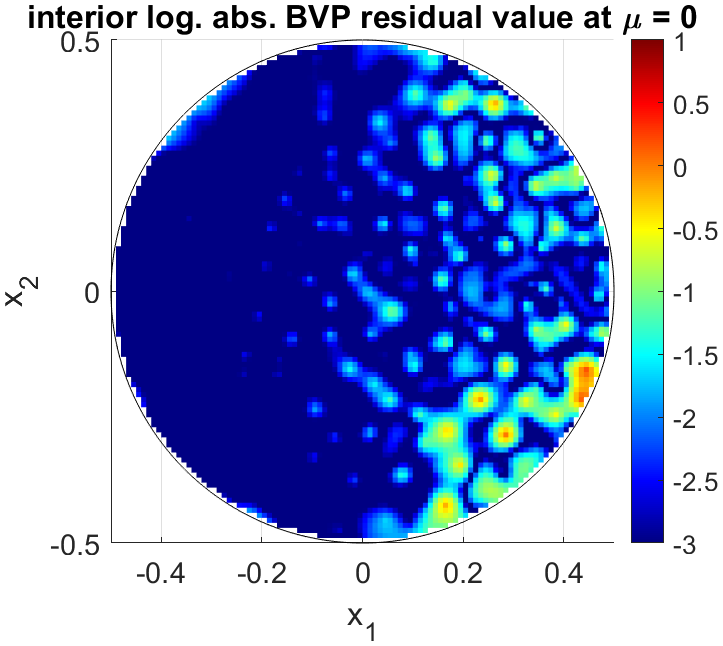}$\quad$
  b) \includegraphics[width=0.50\textwidth]{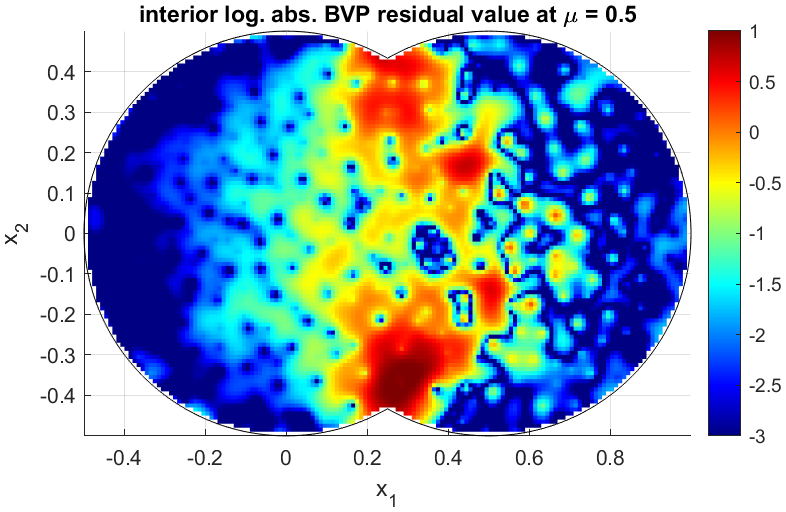} \\
  c) \includegraphics[width=0.7\textwidth]{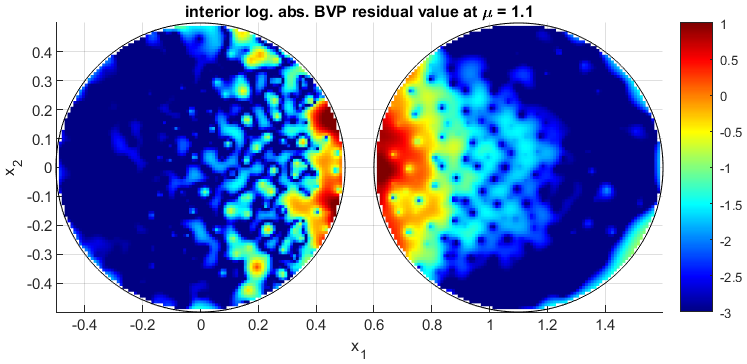}
\end{center}
\caption{\label{fig:moving-circles-nonsmooth-residuals}
  Visualization of interior residuals in case of  the nonsmooth solution for 
  the moving circles example at $\mu=0.0, 0.5$ and $1.1$.
}
\end{figure}

Next we investigate the corresponding BVP functional residuals, by plotting the absolute value of the
interior residuals $|\lambda(u-s_n)|$ for $\lambda\in \Lambda_L$
in Fig.\ \ref{fig:moving-circles-nonsmooth-residuals} for the same parameter values as before.
We observe that there are regions of
considerably
large
absolute residual value. But as well in large regions we observe nicely low residuals.
The problems particularly seem to persist in the regions of the entering corner.
So overall for this nonsmooth example we conclude that
while we do not guarantee any theoretical approximation quality, we are
able to provide some reasonable approximation.

\subsection{Moving Source Term}

\begin{figure}[h!t]
\begin{center}
  a) \includegraphics[width=0.30\textwidth]{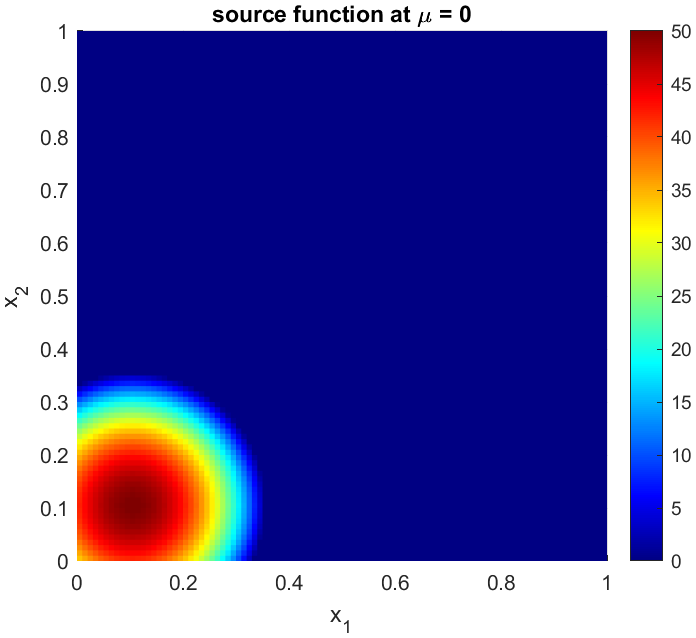}
  \includegraphics[width=0.30\textwidth]{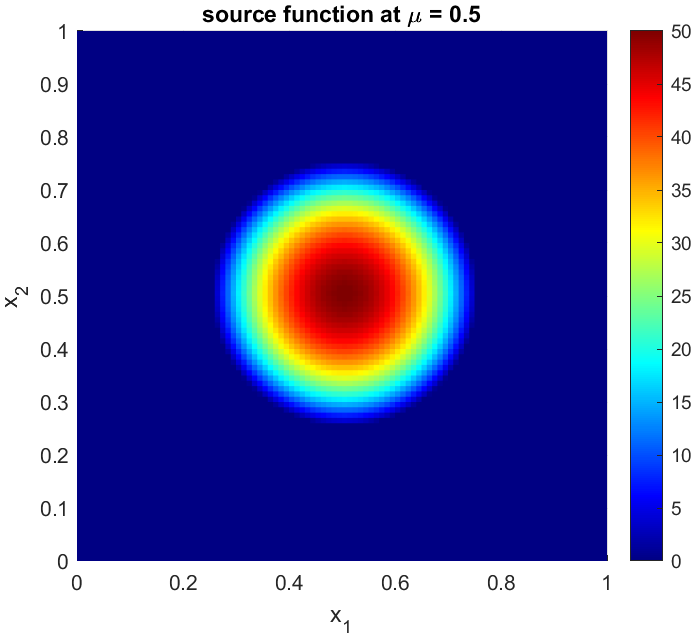}
  \includegraphics[width=0.30\textwidth]{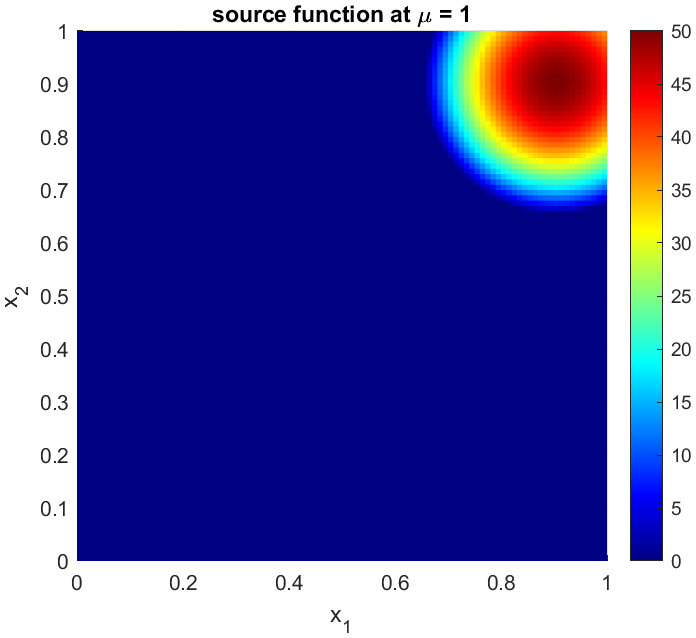} \\
  b) \includegraphics[width=0.30\textwidth]{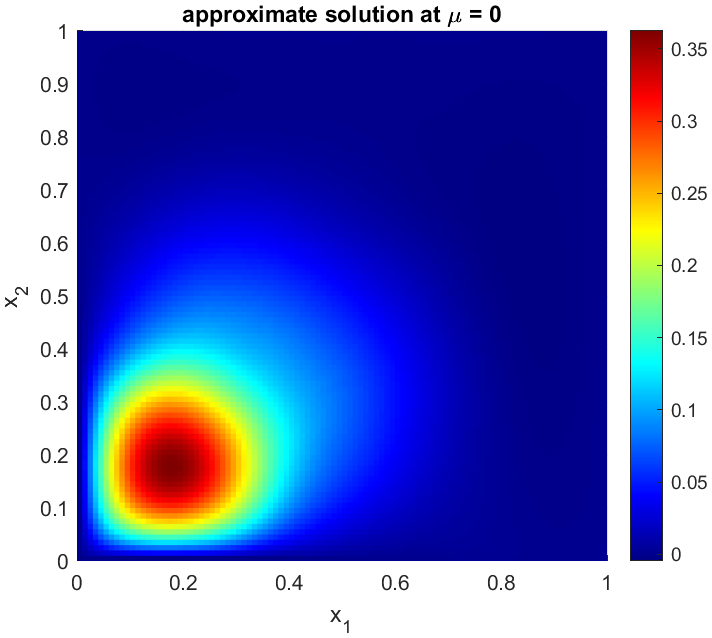}
  \includegraphics[width=0.30\textwidth]{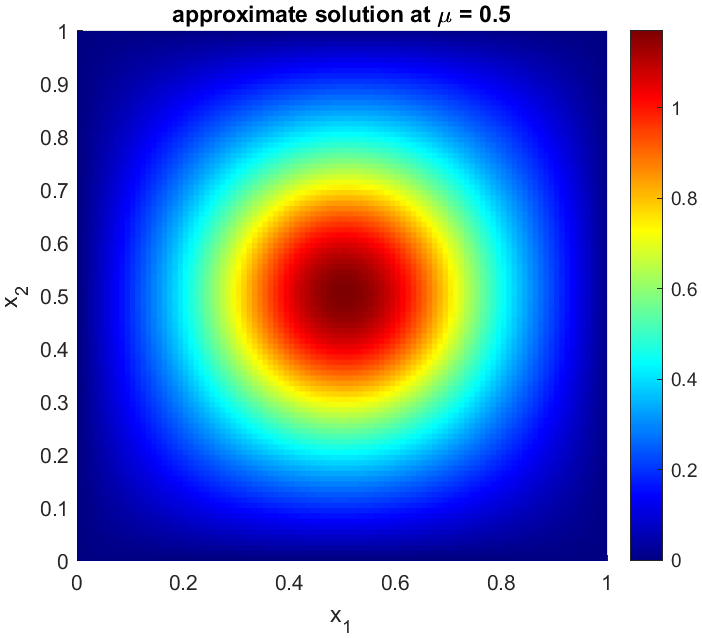}
  \includegraphics[width=0.30\textwidth]{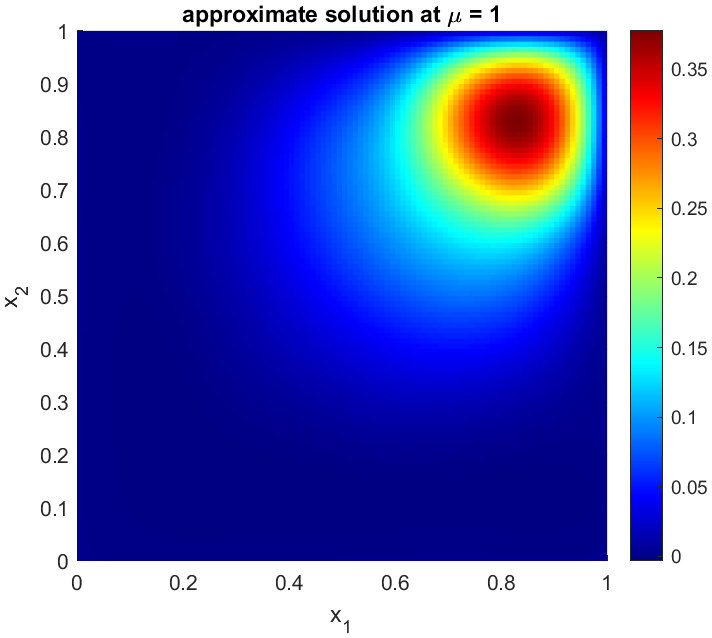}  
\end{center}
\caption{\label{fig:moving-source-and-solution}
  a) Visualization of the moving source functions for parameters $\mu=0.0, 0.5, 1.0$.
  b) Plot of the corresponding approximate solution.}
\end{figure}

The next example is aiming at a moving structure underlying the problem.
Moving structures are typically difficult for RB methods, as this requires
additional (empirical) interpolation steps for the data functions and the
solution manifolds may have slowly decaying Kolmogorov $n$-widths. 
We choose the position domain to be the parameter-independent
unit-square $\Omega(\mu):=\Omega:= (0,1)^2$ for
all $\mu\in \cP:= [0,1]$, hence $d=d_x + d_\mu = 2+1 = 3$.
We define the parametric source function to
be a cut-off inverted parabola $f(x,\mu)= 50 \max(1-\frac{1}{16}\mynorm{x-m(\mu)}^2,0)$
with parametric midpoint $m(\mu):=(\frac{1}{10} + \frac{4}{5} \mu )(1,1)^T \in \Omega$.
This source function is non-affinely parameter dependent, i.e. it is not a
  linear combination of parameter-independent functions with parameter-dependent coefficients.
The source function for the parameters $\mu \in \{ 0.0, 0.5, 1.0 \}$ is visualized in
Fig.\ \ref{fig:moving-source-and-solution}a).
This clearly motivates the
notion ``moving source'' as the source function has a local support, which
is a circle with the given midpoint being translated.
The Dirichlet values are set to homogeneous $g=0$
for all $(x,\mu)\in\partial \bOmega$.

Due to the non-differentiability of $f$ at the boundary of its support, and the
non-differentiable boundary, we cannot expect more than twice differentiability of the solution in the interior
of $\Omega$. Therefore, we again choose the quadratic Mat\'ern kernel for the position and the
parameter coordinate with shape parameters $\varepsilon_x$ and  $\varepsilon_\mu$.
The RKHS of this kernel matches this anticipated regularity in the
interior.  

\begin{figure}[h!t]
  \begin{center}
    a) \includegraphics[width=0.4\textwidth]{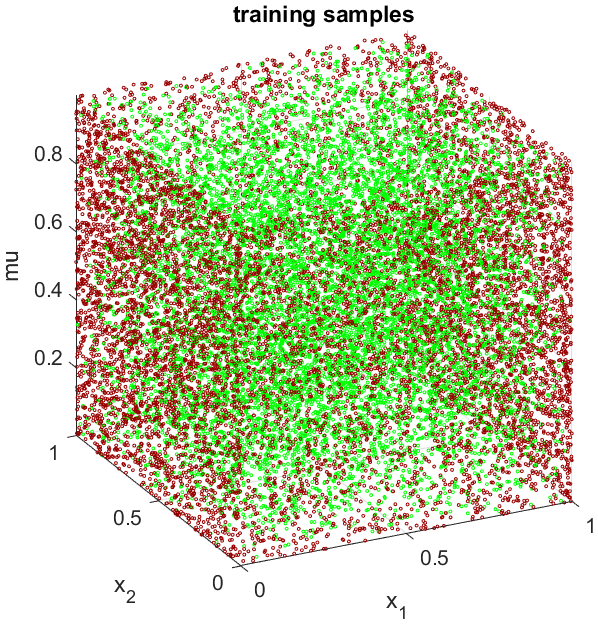} $\quad$
    b) \includegraphics[width=0.4\textwidth]{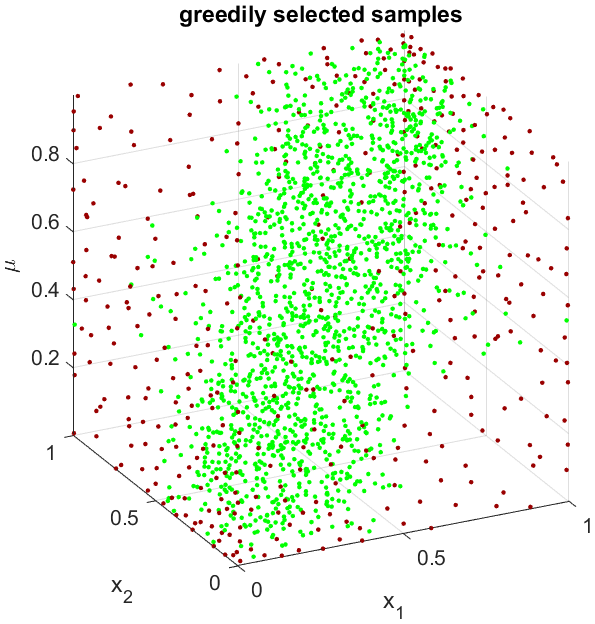} $\quad$ \\
    c) \includegraphics[width=0.4\textwidth]{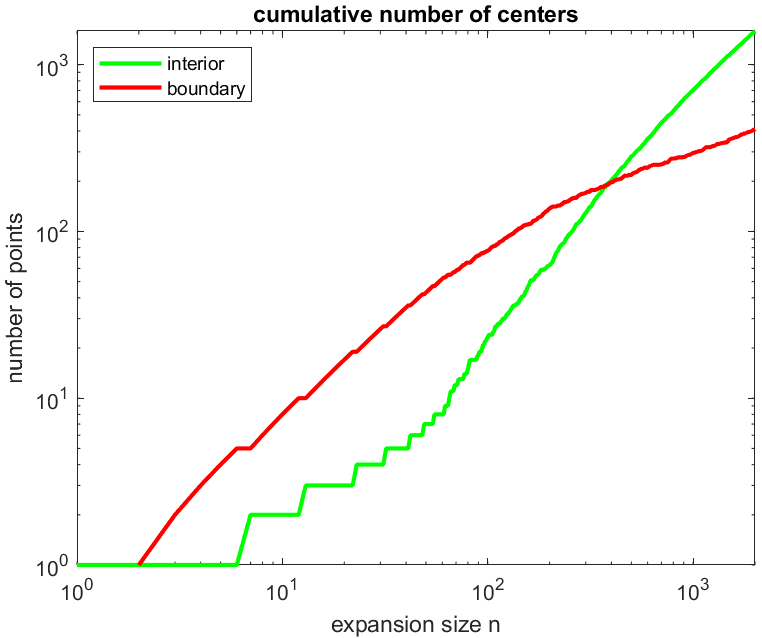}   
    d) \includegraphics[width=0.4\textwidth]{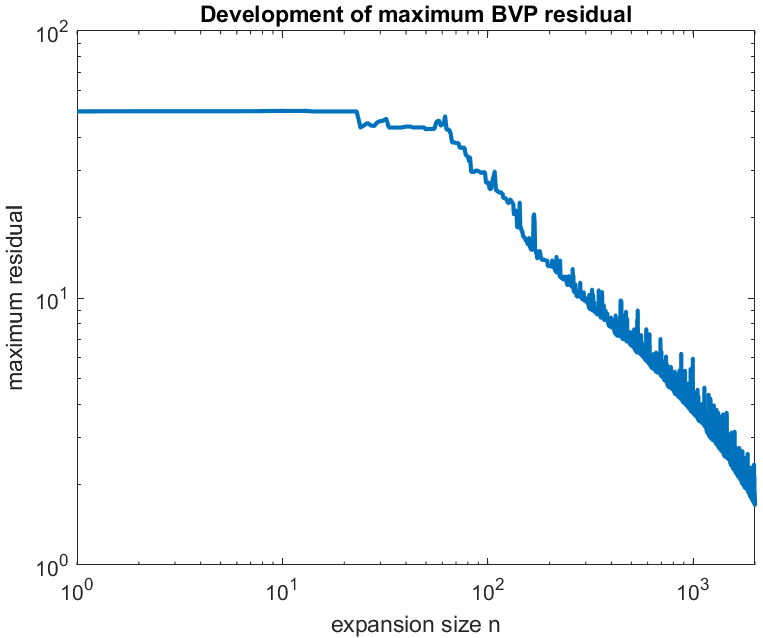}   
\end{center}
  \caption{\label{fig:moving-source-sampling}
    a) Training functionals, indicated by the point evaluation sites
  for the differential (green) and
  boundary operator (red) for 
  the moving source example.
  b) Subset of greedily selected functional samples.
  c) Cumulative number of selected boundary and interior functionals. 
  d) Decay of the maximum training BVP residual absolute value over the iterations. 
}
\end{figure}

Again, we sample randomly and uniformly $N = N_L + N_B =10^4+10^4$ evaluation
points, which correspond to the training functionals as
visualized in Fig.\ \ref{fig:moving-source-sampling}a).
We choose the $f$-greedy strategy with maximum expansion size $n_{max}=2000$ and the same accuracy
and stability tolerances and functional weighting $w_L$ as before.

The two shape parameters and the functional weighting $w_B$ are again
hyper-parameters and were selected by minimizing the loss term of
\eqref{eqn:validation_loss} by 3 consecutive 1D grid-searches.
In this example the default loss weights resulted in a model with too few boundary functional
selections, therefore we reduced the interior weight to $\gamma_L = 0.005$.
The resulting parameters are
$\varepsilon_x = 5.0, \varepsilon_\mu = 15.8114, w_B=1000$.

The approximate solutions for some parameters are visualized in
Fig.\ \ref{fig:moving-source-and-solution}b). 
The symmetry of the problem with respect to the position-diagonal is clearly visible. 
The selected functionals consisting of 1592 interior and 408 boundary functionals are given in Fig.~\ref{fig:moving-source-sampling}b).
Note the clear clustering of the selected functionals around the diagonal
of the cube, which corresponds to the region where the source is nonzero.
The training of the model with these optimized parameters
took 42.03 sec, which indeed is close to a factor 4 increase
compared to the former experiment with an $n=1000$ expansion, hence  
is consistent to the expected $n^2$ scaling. 

The decay of the
maximum BVP residual is plotted in Fig.\ \ref{fig:moving-source-sampling}d), which
again indicates an overall rate which seems algebraic, but again
shows clear non-monotonicity.
The cumulative number of selected interior and boundary functionals is shown in
Fig.\ \ref{fig:moving-source-sampling}c).
Again, we observe different phases: The initial plateau of up to $n=50$ is caused by the
choice of mostly boundary functionals due to the large value of $w_B$.
The BVP functional error only 
decays adequately as soon as a larger number of interior
functionals have been chosen. 

  In order to obtain more quantitative error statements, we evaluate
  the $L^\infty(\Omega_x)$-error 
  compared to a reference solution obtained by a FEM.
  For this we choose a uniform triangular mesh 
  with $20000$ elements and P1 ansatz functions.
  In addition to the PDE-$f$-greedy procedure ($\beta=1$) we now also
  include the $\beta=0.75, 0.5$ and $\beta=0$ cases.
  Also, we consider a full kernel collocation and a coarser FEM approximation as methods for comparison.
  For a random position-parameter test set, all parameter components in the test vectors would differ and require
  a separate fine FEM solve for each of those parameter values, which would clearly be computationally infeasible.
  Instead, we construct a structured Cartesian product test set based on choosing
  $n_\mu=20$ equidistant parameter values $\cP_{h} \subset \cP$ and $n_x=101^2$ position samples as the
  vertices $\Omega_h \subset \Omega$ of the fine reference FEM mesh. 
  Table \ref{tab:table-moving-source}
  lists the resulting parameters and errors.

  \begin{table}[htp!]
    \begin{tabular}{|c||c|c|c|c|c|c|}
      \hline
      \multicolumn{7}{|c|}{PDE-$\beta$-greedy} \\
            \hline
      $\beta$ & $n$ & $\varepsilon_x $ &   $\varepsilon_\mu$ &  $w_B$ & $r_{bnd}$  & $L^\infty$-error \\
      \hline 
   $0$  & $2000$ & $5$ & $5$ & $100$ & $0.2825$ & $0.17416$ \\
   $0.5$  & $2000$ & $5$ & $15.8114$ & $200$ & $0.2795$ & $0.032093$ \\ 
   $0.75$  & $2000$ & $5$ & $15.8114$ & $100$ & $0.1135$ & $0.027812$ \\ 
    $1$  & $2000$ & $5$ & $15.8114$ & $1000$ & $0.204$ & $0.032785$ \\
      \hline 
      \hline
      \multicolumn{7}{|c|}{Full collocation} \\
            \hline
       & $n$    & $\varepsilon_x $ &   $\varepsilon_\mu$ &     & $r_{bnd}$  & $L^\infty$-error \\
      \hline 
       & $2000$ & $5$              &    $5$             &      & $0.1$    & $0.21286$ \\
      \hline 
      \hline      
      \multicolumn{7}{|c|}{FEM} \\
            \hline
      $n_x$ & $n=n_\mu \cdot n_x$ & \multicolumn{4}{c|}{} & $L^\infty$-error \\
      \hline
      $100$ & $2000$ & \multicolumn{4}{c|}{}
      & $0.088282$ \\ 
      \hline
      \end{tabular}
    \caption{Resulting selected parameters and $L^\infty$-errors for different
      PDE-$\beta$-greedy variants: Full collocation and FEM for the moving source example.
      \label{tab:table-moving-source}}
  \end{table}
  
  To make the FEM and PDE-$\beta$-greedy approaches comparable we choose models with
  similar number of degrees of freedom (DOFs),
  where the number of DOFs for a nonparametric FEM solution are simply multiplied by the
  number of parameter samples, as this is the amount of coefficients that must be stored to represent
  the $n_\mu$ parametric solutions.
  Therefor we consider an equidistant $10\times 10$ nodes FEM mesh on the unit square and
  choose the $20$ parameter samples, as this results in 2000 degrees of freedom for those 20 FEM solutions.
  For the full kernel collocation approach, only the kernel parameters were selected by validation,
  the ratio of boundary value points was fixed a priori.
  We can observe that all the greedy procedures with $\beta>0$ are clearly superior over
  the full collocation approach, which is not
  surprising as the latter uses a fixed non-adapted center set, while they are better but
  comparable for $\beta=0$, which still selects the centers from a larger set.
  In comparison to the FEM approach, we can recognize that the target function-dependent greedy selection criteria result
  in models with clearly lower error than the FEM approximation with comparable
  number of DOFs, even if it is not a large factor of improvement. So the function adaptivity pays off.
  Note that the FEM error is to be understood as being optimistic for the FEM
  method, as for a full parametric assessment, some parametric interpolation, etc.
  would need to be included, which cannot decrease the FEM errors for the snapshot solutions but will presumably yield
  larger values for intermediate parameter values.
  Especially this FEM error is a lower bound for any parametric MOR method that is based
  on those 20 parametric FEM solution snapshots and has a reproduction property, e.g.\ 
  parametric interpolation, Galerkin-projection, etc.

\subsection{High-dimensional Domains}

To investigate the possibilities of the scheme for high-dimensional cases we consider as positional domain
the non-parametric unit-hypercube $\Omega:=(0,1)^{d_x}$, $d_x \in \{2,\ldots,9\}$ with single
scalar parameter $\mu\in [0,1]$, hence $d_\mu=1$ and overall dimension $d\in\{3,\ldots,10\}.$
We consider a model that can easily be generalized to arbitrary input dimension by
prescribing as known solution a sinus wave 
$u(x,\mu):=\sin(\scalprod{x}{\kappa(\mu)})$
with parametrized frequency
$$\kappa(\mu)=\frac{1}{d_x}((1-\mu)\kappa_0 + \mu \kappa_1)(1,\ldots,1)^T \in \R^{d_x},$$
with $\kappa_0:=\pi$ and $\kappa_1:=2\pi$ such that $\mu=1$ always results in a complete
sinus period and $\mu=0$ in half a period along the diagonal of the hypercube independently of $d_x$.
The RHS functions $f$ and $g$ are chosen correspondingly. 
This again allows us to assess the true pointwise errors of the approximation to the exact solution.
To cover higher dimensions more accurately we raise the number of training
samples of functionals to 
$N_L:=N_B:=10^5$, 
hence $N=2\cdot 10^5$. 

In view of the $C^\infty$ regularity of the solution, as before we use the
product of a position and parameter Gaussian kernel.

We use an $f$-greedy criterion and fix the maximum expansion size to $n_{max} = 1000$ and the stopping tolerances
as before. The functional weightings were set to $w_L:=1$ and $w_B:=100$.

For comparability over the dimensions, we do not do a meta-search of parameters. Instead we
set $\varepsilon_x^2:=5.0/d_x, \varepsilon_\mu^2:= 2.5$ 
which yields reasonable but
certainly improvable results in all considered dimensions.
As the diagonal of the unit-cube $\Omega$ scales with $\sqrt{d_x}$, the specified
  choice of $\varepsilon_x$ results in ``equal'' kernel decay along the
  diagonal in all dimensions. We 
refer to this as ``isotropic Gaussian'' in the following. 
Note that a sufficiently fast decaying choice of dimension-dependent shape
parameters is also the key for obtaining dimension-independent constants in
worst-case error bounds \cite{Fasshauer2012c}. 

\begin{figure}[h!t]
\begin{center}
  a) \includegraphics[width=0.45\textwidth]{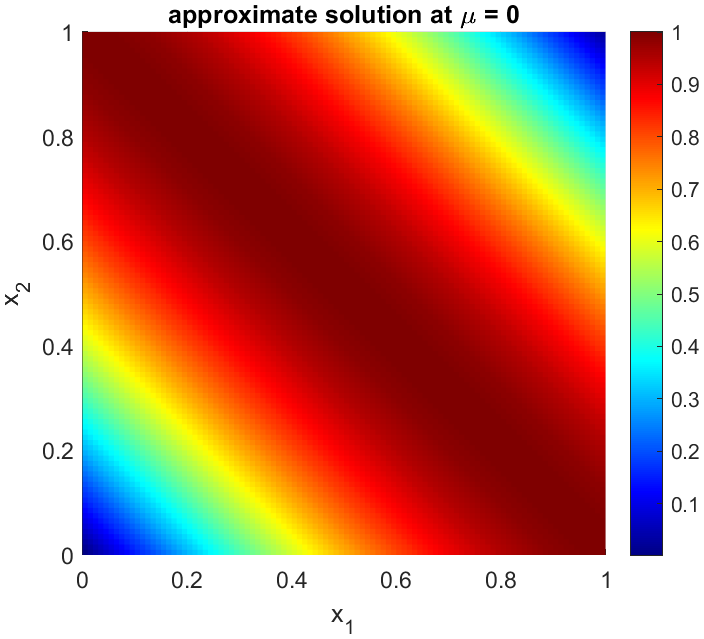}
  b) \includegraphics[width=0.45\textwidth]{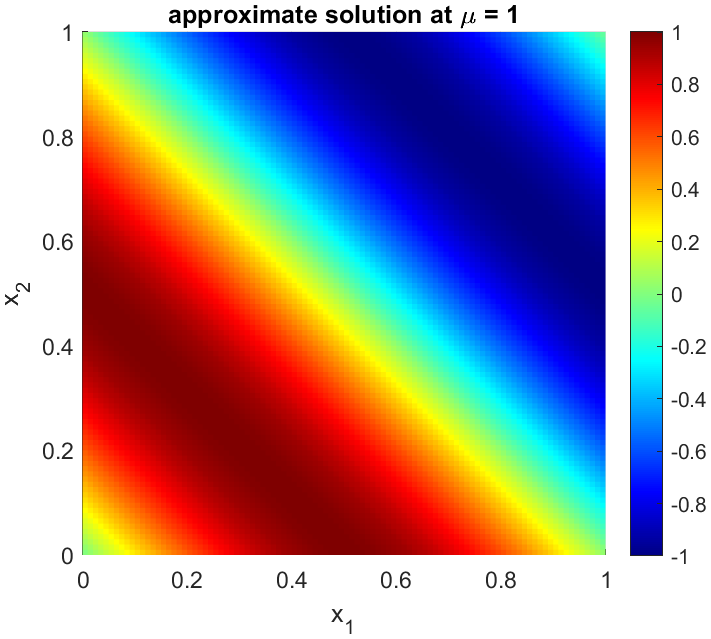} 
\end{center}
\caption{\label{fig:nd-sinus-source-solution}
  Visualization of sinus-source example approximate solution in $\R^{d_x}$ for $d_x=2$ a) at
  $\mu=0$ and b) at $\mu=1$.}
\end{figure}

In Fig.\ \ref{fig:nd-sinus-source-solution} we illustrate the parametric approximate
solution evaluated for two extreme parameters a) $\mu=0$ and b) $\mu=1$ for the $d_x=2$ case.
The doubling in frequency along the diagonal is clearly visible. 
The accurate capturing of the diagonal symmetry even with an isotropic kernel can be
verified.  

\begin{figure}[h!t]
\begin{center}
  a) \includegraphics[width=0.45\textwidth]{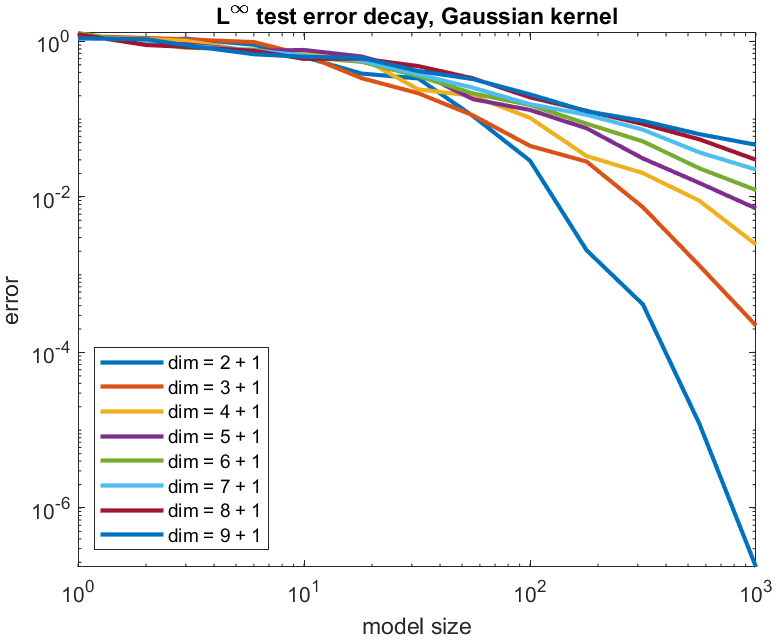} 
  b) \includegraphics[width=0.45\textwidth]{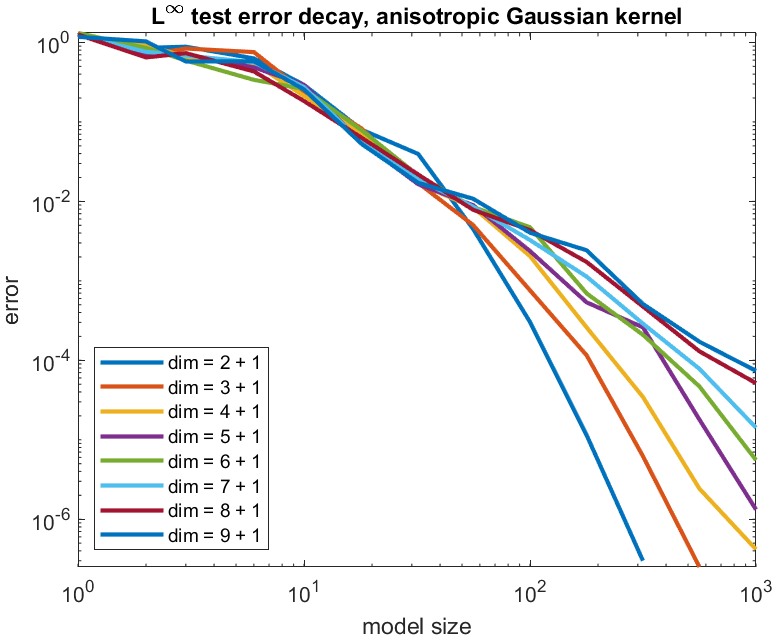} 
\end{center}
\caption{\label{fig:nd-sinus-source-convergence}Test 
  $L^\infty$-error convergence for the high-dimensional sinus-source example with varying position-dimensions and using
  a) an isotropic Gaussian kernel and b) an anisotropic Gaussian kernel.}
\end{figure}

In order to illustrate the $L^\infty$ convergence rates, we evaluated the error measured over
$n_{L,test} + n_{B,test} = 10000+10000$ test functionals for increasing model sizes sampled at 13 (approximately) logarithmically
equally distant model sizes $n=10, \ldots, 1000$. The resulting errors for the different dimensions
are plotted in Fig.\ \ref{fig:nd-sinus-source-convergence}a).
For every dimension we clearly recognize an exponential decay in the log-log plot.
With increasing space dimension it is more difficult to provide
small approximation errors, which is reflected by the slower decay.

In the given example, the radial kernel obviously is not adapted to
the solution structure. As final experiment, we want to illustrate how problem-adapted
kernels can lead to considerable improvements.

Using strong prior knowledge about the translation-invariance of the
solution, we want to reflect this with an anisotropic kernel choice for the position variable.
For this, we choose a suitable matrix $B$
in \eqref{eqn:anisotropic-gaussian} by splitting the position
space into $v_1:= \frac{1}{\sqrt{d_x}}(1,\ldots,1)^T \in \R^{d_x}$ and an orthogonal complement,
i.e.\ $v_2,\ldots v_{d_x} \perp v_1$ being pairwise orthogonal and normalized $\mynorm{v_i}=1$. 
Then we obtain the orthogonal matrix $V:=[v_1,\ldots,v_{d_x}] \in \R^{d_x \times d_x}$ and set
\begin{align*}
  B:= V \left(\begin{array}{cc} 5.0/d_x & 0 \\ 0 & 0.0025\cdot I_{(d_x-1) \times (d_x-1)}\end{array}\right) V^T.
\end{align*}
The parameter kernel is chosen as before as standard Gaussian with $\varepsilon_\mu^2 =2.5$.
  
The resulting $L^\infty$-error convergence curves are
given in Fig.\ \ref{fig:nd-sinus-source-convergence}b).
Compared to plot Fig.\ \ref{fig:nd-sinus-source-convergence}a), we indeed see a considerable improvement
in approximation, which even seems to be partially dimension-independent
until certain model sizes are reached.  

We want to comment that the assumption of knowing the most important direction
of variation, i.e.\ $v_1$ can be relaxed. 
By applying an active subspace detection procedure \cite{Co2015}, i.e.\ 
sampling the source function gradient at the training points and computing
$v_1$ as the principal component, an active subpace for $f$ can be identified.
When using this numerically identified subspace $v_1$ in the anisotropic kernel,
the resulting $L^\infty$-errors only change marginally in the present example.
A reason for the success of this procedure may be that in our case
due to the isotropy of the differential operator the active subspace 
from the source function transfers to the active subspace of the solution.

These results -- even if not being meta-parameter-optimized -- indicate promising approximation properties
for high-dimensional examples.

\section{Summary and Conclusion}

\label{sec:conclusion}
In the present article we extend the scale of $\beta$-greedy schemes for PDE approximation
to the parametric case by the use of position-parameter kernels.
The procedure is easily implemented as a greedy kernel interpolation scheme using
an abstract kernel on BVP evaluation functionals.
Correspondingly the method inherits computational advantages such as
linear training complexity in the overall number of candidate collocation
centers and only quadratic complexity in the final kernel expansion size. 
Analytical convergence statements suggest that one can expect exponential
convergence rates, hence very efficient and sparse expansions, in the case of
smooth solutions and smooth kernel. Especially for $\beta>0$ an additional
dimension-independent decay factor appears suggesting that
RHS-adaptive center selection is indeed provably preferable over
selection schemes that are agnostic of the special source/boundary values.  
The exponential convergence rates were
experimentally verified in different scenarios.
We do not suggest the present approach as a competitor to standard RB methods
for parameter separable low space-dimensional elliptic cases, where
RB-methods are known to provide exponential convergence rates.
Rather we demonstrated that the present PDE-$\beta$-greedy approach
can be beneficial in cases that
cover nontrivial geometry parametrization including topology changes,
moving structures or high-dimensional state domains.
We do not require explicit parametric geometry transformations of the
domain as long as we can sample points from the interior and the boundary
to generate the candidate collocation point sets. 

Beside those positive insights, limitations of the procedure appear in
cases, e.g.\ with singular solutions. This is due to the fact that we
aim for strong solutions, hence assume classical regularity.
Problems with solutions that have corner singularities, cannot be approximated
well with the presented schemes. Still
some limited approximation is possible, cf.\ the singular solution case of Sec.\ \ref{sec:moving-circles}
and the numerical comparison to FEM
approximation in \cite{PDEgreedy2025}.

As interesting outlook we see experimental investigation for
more general non-elliptic PDEs, e.g.\ advection or wave-equation problems.
These also pose problems for linear subspace-based
MOR, as they face slowly decaying Kolmogorov $n$-widths \cite{greif2019decay,ohlberger2016reduced}
when relying on separating the spatial and time coordinate and
aiming for a spatial approximation space that is used for the whole time axis.
In contrast, the current kernel-greedy schemes would circumvent this problem by inherently
constructing efficient approximation spaces in the time-space domain.
This will be addressed in more detail in future work. 

A further immediate direction for extension is the treatment of 
nonlinear PDEs by iterative linearizations and corresponding
fixpoint schemes. This would then enable efficient solution of
problems in high-dimensions such as kinetic equations, optimal control
(Hamilton-Jacobi-Bellmann equatios), finance (option pricing) or high-dimensional
parameter spaces.

Analytically, it may be interesting to investigate the case of ``escaping the native space''~\cite{Narcowich2006o},
i.e.\ the case where the target function is coarser than the used kernel
and it may be interesting to see how the convergence rates degrade depending on
the difference in the smoothness between the kernel and the solution. 

Methodologically, it seems evident that problem specific kernels
have the potential to provide very compact models.
Thus it seems immediate to extend the present greedy kernel collocation approach
with multilayer or deep-kernel learning in the spirit of \cite{WHKOS24,WMP24}. 
For problems with structural invariances, it may be promising to
use structure-preserving kernels in the position component, e.g.,
divergence-free kernels~\cite{Fuselier2009,Wendland2009} for incompressibility conditions. 

\section*{Acknowledgements}
The authors acknowledge the funding of the project by the Deutsche
Forschungsgemeinschaft (DFG, German Research Foundation) under number 540080351 and
Germany's Excellence Strategy - EXC 2075 - 390740016.

\noindent

\begin{appendices}

\section{Proofs of Propositions}\label{secA1}

\begin{delayedproof}{prop:algebraic-rates-functional-residuals}
  Without loss of generality, we assume
    \begin{equation}
      P_{\Lambda_0}(\lambda)\leq 1 \mbox{ for all } \lambda\in \Lambda. \label{eqn:PLambda0leq1}
    \end{equation}
  If this condition is not met, the boundedness of $\Lambda$ implies the existence of a
  factor such that the scaling of the functionals ensures
  \eqref{eqn:PLambda0leq1}.
  The proof then continues with the scaled functionals
  and the scaling factor enters the final constant $C$ in the decay rate.  
  With \eqref{eqn:PLambda0leq1} the assumptions of Thm.\ 17 in
  \cite{PDEgreedy2025} are valid and 
  from its proof we cite
  the bound 
    (using $\bar n$ and $\bar C$ to discriminate from the
  variables in the proposition) 
  \begin{equation}
  \min_{i=\bar n+1,\ldots,2\bar n} \sup_{\lambda \in \Lambda} |\lambda(u-s_i)|
  \leq \bar C {\bar n}^{-\beta/2} {\bar n}^{-\alpha} \mynorm{e_{\bar n+1}}_{\nsbar} \log(\bar n)^\alpha,
  \label{eqn:beta-greedy-algebraic}
  \end{equation}
   which is valid for $\bar n\geq 3$. 
  As $e_{\bar n+1} = u-s_{\bar n+1}$  and $s_{\bar n+1}$ is
    obtained as an orthogonal projection we
    conclude by Pythagoras that $\mynorm{e_{\bar n+1}}\leq \mynorm{u}$.
    On the left hand side the minimum can be lower bounded by taking
    the minimum over all $i=1,\ldots,2\bar n$ instead of considering the
    smaller range $i=\bar n+1,\ldots 2\bar n$.
    We abbreviate
        $\gamma := -\beta/2 + 1/2 -(\tau-2)/d < 0$.
    For odd $n=2\bar n + 1$ we obtain
    \begin{equation}
      \bar n^\gamma= \left(\frac {n-1}{2}\right)^\gamma 
      = \left(\frac12\right)^\gamma \left(1-\frac 1n\right)^\gamma n^\gamma.
      \label{eqn:aux4}
    \end{equation}
    Assuming $\bar n\geq 3$ we have $n\geq 7$ and can upper bound the middle
    factor by $\left(\frac{6}{7}\right)^\gamma$ and define the new constant as
    $$C:=\bar C \left(\frac12\right)^\gamma \left(\frac{6}{7}\right)^\gamma=\bar C \left(\frac{3}{7}\right)^\gamma.$$
    Further the logarithm factor
    can be upper bounded by $\log(\bar n) < \log(2\bar n+1) = \log(n)$.
    Thus we overall verify the statement for odd $n$ as
    $\min_{i=1,\ldots, n}(\cdot) = \min_{i=1,\ldots,2\bar n + 1} (\cdot) \leq \min_{i=1,\ldots,2\bar n} (\cdot)$
    and using \eqref{eqn:beta-greedy-algebraic} with \eqref{eqn:aux4}.
    For even $n=2\bar n$ we similary obtain 
    \begin{equation}
      \bar n^\gamma= \left(\frac n2\right)^\gamma = 
      \left(\frac n2\right)^\gamma= \left(\frac12\right)^\gamma n^\gamma
      \label{eqn:aux3}
      \end{equation}
      which then is similarly resembled by the new constant $C$ and overall gives the statement
      for all $n\geq 6$.      
\end{delayedproof}

\begin{delayedproof}{prop:exponential-rates-functional-residuals}
  We argue very analogously to the previous proof, hence shorten the presentation: 
  We start with a bound obtained in \cite{vogel24} using $\bar n$, etc.\ to discriminate from the present notation:
  $$
  \min_{i=\bar n + 1,\ldots, 2 \bar n} \sup_{\lambda \in \Lambda} |\lambda(u-s_i)| \leq \bar C \bar n^{-1/2}
      \mynorm{u-s_{\bar n+1}}_{\nsbar} e^{-\bar c_1 \bar n^{\frac{1}{2d}}}.
      $$
  We can lower bound the RHS by extending the index range from $i=1,\ldots, 2 \bar n$.
  The error norm on the RHS can again be upper bounded by $\mynorm{u}_{\nsbar}$.
      For the even case $n:=2\bar n$ the exponent in the last factor on the RHS can be
      rewritten as
      $$
         -\bar c_1 \bar n^{\frac{1}{2d}} = 
         -\bar c_1 \left(\frac n2 \right)^{\frac{1}{2d}} = -\bar c_1 \left(\frac 12 \right)^{\frac{1}{2d}}
         n^{\frac{1}{2d}}
         = -c_1 n^{\frac{1}{2d}}
         $$
      by setting $c_1:= \bar c_1 \left(\frac 12 \right)^{\frac{1}{2d}}$.
      As in \eqref{eqn:aux3} but using a different $\gamma$
      the factor $\bar n^\gamma$ can be expressed by
      a corresponding constant factor and $n^\gamma$.
      In the case of odd $n=2\bar n + 1$ similar reformulations and upper bounding
      according to \eqref{eqn:aux4} can be applied requiring that $n$ is sufficiently large,
      e.g.\ $n\geq 6$.  
      The statement then follows by suitable definition of the constant $C$. 
\end{delayedproof}

\end{appendices}

\end{document}